\documentclass{amsart}

\newcommand{\PSbox}[3]{\mbox{\rule{0in}{#3}\includegraphics{#1}\hspace{#2}}}
\theoremstyle{plain}
\newtheorem{thm}{Theorem}
\newtheorem{cor}[thm]{Corollary}
\newtheorem{lem}[thm]{Lemma}
\newtheorem{prop}[thm]{Proposition}

\numberwithin{equation}{section}

\newcommand{\thmref}[1]{Theorem~\ref{#1}}
\newcommand{\propref}[1]{Proposition~\ref{#1}}
\newcommand{\corref}[1]{Corollary~\ref{#1}}

\newcommand{\lemref}[1]{Lemma~\ref{#1}}

\newcommand{\Z}    {{\mathbb Z}}
\newcommand{\C}    {{\mathbb C}}

\newcommand{\Cr}   {{\textup{Cr}}}
\newcommand{\Pl}   {{\textup{P}}}
\newcommand{\Plt}  {{\textup{\~P}}}

\newcommand{\Exp}  {{\textup{Exp}}}
\newcommand{\Var}  {{\textup{Var}}}
\newcommand{\Prob} {{\textup{Prob}}}

\title[Local statistics for random domino tilings]{Local statistics
for random domino tilings of the Aztec diamond} 
\author{Henry Cohn}
\address{Department of Mathematics, Harvard University}
\email{cohn@math.harvard.edu} 
\author{Noam Elkies} 
\address{Department of Mathematics, Harvard University} 
\email{elkies@math.harvard.edu}
\author{James Propp} 
\address{Department of Mathematics, MIT}
\email{propp@math.mit.edu} 
\thanks{The second author was supported in part by an NSF Presidential
Young Investigator Award and by a Fellowship from the Packard
Foundation. The third author was supported in part by NSA grant
MDA904-92-H-3060 and NSF grant DMS 9206374, and by a career
development grant from the M.I.T.\ Class of 1922.} 
\subjclass{Primary 60K35, 82B20; Secondary 05A16, 60C05}
\date{April, 1996}
\begin{document}

\begin{abstract}
We prove an asymptotic formula for the probability that, if one
chooses a domino tiling of a large Aztec diamond at random according
to the uniform distribution on such tilings, the tiling will contain a
domino covering a given pair of adjacent lattice squares.  This
formula quantifies the effect of the diamond's boundary conditions on
the behavior of typical tilings; in addition, it yields a new proof of
the arctic circle theorem of Jockusch, Propp, and Shor.  Our approach
is to use the saddle point method to estimate certain weighted sums of
squares of Krawtchouk polynomials (whose relevance to domino tilings
is demonstrated elsewhere), and to combine these estimates with some
exponential sum bounds to deduce our final result.  This approach
generalizes straightforwardly to the case in which the probability
distribution on the set of tilings incorporates bias favoring
horizontal over vertical tiles or vice versa.  We also prove a fairly
general large deviation estimate for domino tilings of
simply-connected planar regions that implies that some of our results
on Aztec diamonds apply to many other similar regions as well.
\end{abstract}

\maketitle

\section{Introduction}
\label{sec-introduction}

\subsection{Statement of the main theorem.}
\label{ssec-statement}

Random domino tilings of finite regions often exhibit surprising
statistical heterogeneity.  Such heterogeneity would be expected in
the vicinity of the boundary, but in fact the presence of a boundary
can make its influence felt well into the interior of the region.  The
research that led to this article is part of an ongoing effort to
understand this phenomenon.  The results proved here are the first to
give a precise description of how local statistics for domino tilings
can vary continuously throughout a region in response to the
imposition of specific boundary conditions.

Those who study random tilings of finite regions (in the plane) by
dominos have tended to focus on regions that are rectangles of even
area.  In particular, Burton and Pemantle \cite{local} have done an
intensive analysis of the small-scale structure of such tilings.
Their work shows that once one gets away from the boundary of the
rectangle, random tilings tend to exhibit statistical isotropy.  Among
all random processes that take their values in the set of domino
tilings of the plane, the Burton-Pemantle process has maximal entropy,
and it is unique in this regard; for this reason alone, it is worth
further study.

However, if one looks at random domino tilings of tileable finite
regions in general, one finds that local behavior far from the
boundary need not be governed by maximal entropy statistics, but can
look very different.  Moreover, the local behavior seen in one part of
the region is in general different from local behaviors seen
elsewhere.

\begin{figure}
\begin{center}
\PSbox{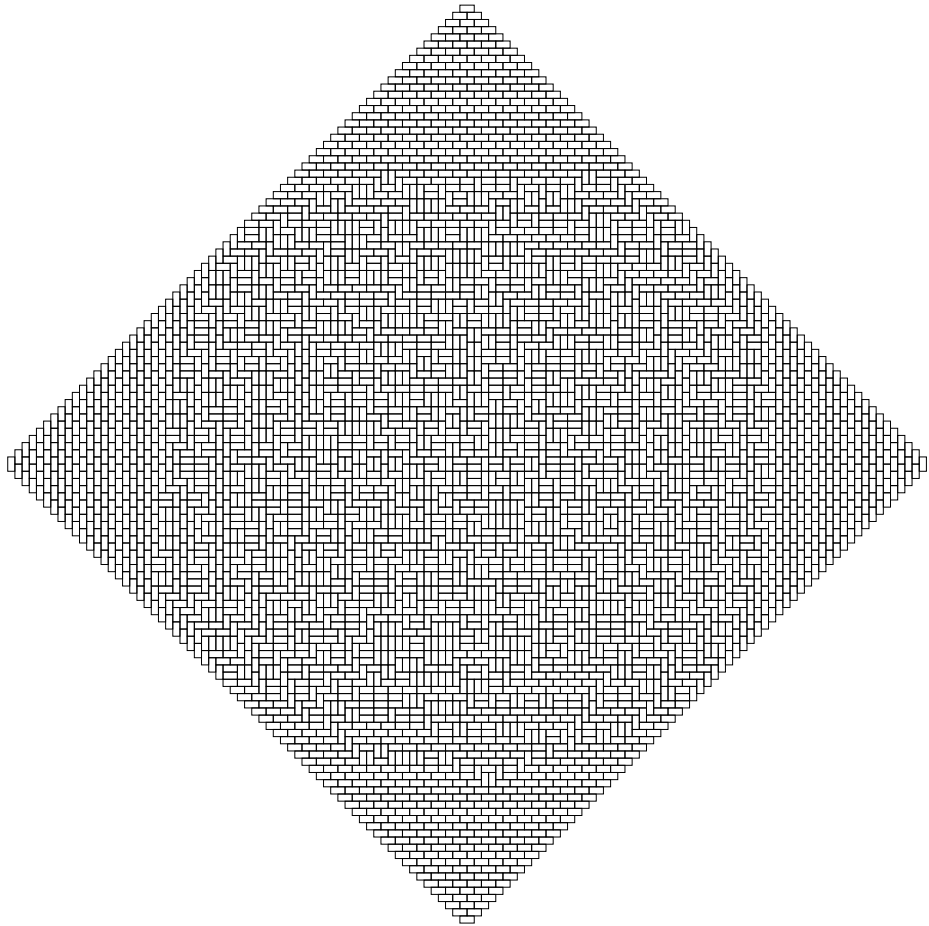 hoffset=0 voffset=0}{350pt}{350pt}
\end{center}
\caption{A random domino tiling of an Aztec diamond of order 64.}
\label{fig-aztec}
\end{figure}

One especially tractable proving ground for the study of this
statistical heterogeneity has been the family of finite regions known
as Aztec diamonds, introduced and studied in \cite{alternating}.
Figure~\ref{fig-aztec} shows an Aztec diamond of order 64 tiled
randomly by dominos.  In general, the Aztec diamond of order $n$ can
be defined as the union of those lattice squares whose interiors lie
inside the region $\{(x,y): x+y \leq n+1\}$.

It was shown in \cite{circle} (and will be proved in
subsection~\ref{ssec-ACT}~by different methods) that, asymptotically,
the circle inscribed in the Aztec diamond of order $n$ serves as a
boundary between domains of qualitatively different behavior.  We call
this circle the arctic circle, because, as one can see from
Figure~\ref{fig-aztec}, the dominos outside the arctic circle are
frozen into a brickwork pattern.  To state the theorem more precisely,
we impose a checkerboard coloring on the Aztec diamond of order $n$,
so that the leftmost square in each row in the top half of the diamond
is white.  We say a horizontal domino is {\sl north-going} or {\sl
south-going} according to whether its leftmost square is white or
black, and we say a vertical domino is {\sl west-going} or {\sl
east-going} according to whether its upper square is white or black.
(The motivation for this terminology comes from the ``domino
shuffling'' algorithm introduced in \cite{alternating} and used in
both \cite{circle} and \cite{gip}; this algorithm permits one to
generate random domino tilings of Aztec diamonds in such a way that
every possible tiling has the same probability of arising as every
other, and indeed it was this algorithm that we used to generate the
tiling shown in Figure~\ref{fig-aztec}.)

Say that two dominos are adjacent if they share an edge (i.e., their
boundaries overlap on a segment of length $1$ or more), and say that a
domino is adjacent to the boundary of the Aztec diamond if it shares
an edge with the boundary.  We define the {\sl north polar region} as
the union of those north-going dominos that are each connected to the
boundary by a sequence of adjacent north-going dominos.  The south,
west, and east polar regions are defined similarly, and the {\sl
temperate zone} is the union of those dominos that belong to none of
the four polar regions.

The arctic circle theorem of \cite{circle} states that for every
$\varepsilon > 0$, if one takes $n$ sufficiently large, then for all
but an $\varepsilon$ fraction of the domino tilings of the diamond of
order $n$, the border of the temperate zone stays within distance
$\varepsilon n$ of the circle of radius $n/\sqrt{2}$ with center
$(0,0)$.  In particular, this implies that if one increases the radius
of the circle by $\varepsilon n$, then with probability greater than
$1-\varepsilon$, in each of the four regions in the Aztec diamond that
lie outside the enlarged disk, all dominos are aligned with their
neighbors in brickwork patterns.  The theorem also implies that if one
decreases the radius of the disk by $\varepsilon n$, then within the
shrunken disk dominos with different orientations are in some sense
interspersed among one another (with probability greater than
$1-\varepsilon$); however, the theorem by itself gives no information
on their distribution.

In \thmref{main} of this article we will give a quantitative analysis
of the behavior of random tilings in the inner, disorderly zone.  In
particular, we will give an asymptotic formula for the proportion of
domino tilings of the Aztec diamond of order $n$ that contain a domino
at a specified location, i.e., the {\sl placement probability} for
that location.  This formula depends only on the orientation of the
domino, its parity relative to the natural checkerboard coloring of
the Aztec diamond, and the relative position of the domino within the
Aztec diamond (in normalized coordinates).  One consequence of our
formula is that random domino tilings exhibit ``total statistical
heterogeneity'' within the central zone.  That is to say, any two
patches within the temperate zone that are macroscopically separated
(i.e., separated by a distance on the order of $n$) will exhibit
distinct statistics.  (For a precise statement, see
subsection~\ref{ssec-heterogeneity}.)

Our work builds on the generating functions derived in \cite{gip}.
One of them is a rational function in three variables whose
coefficients are the placement probabilities for which an asymptotic
formula is sought.  The authors of the earlier article carried out a
relatively straightforward complex integration to calculate
coefficients corresponding to dominos in the $2\times2$ block in the
middle of the Aztec diamond; the resulting exact formula implies that
in a diamond of order $n$, these placement probabilities are $\frac14
+ O(\frac1n)$.  In the present article we will apply the saddle point
method to estimate contour integrals associated with more general
coefficients of a related generating function (also derived in
\cite{gip}).

We can now prepare to state our main result.  We call the union of two
adjacent squares in the Aztec diamond a {\sl domino space\/}, to avoid
confusion between actual dominos occurring in a particular tiling and
the locations in which dominos can occur.  Domino spaces are
classified as north-going, south-going, west-going, or east-going in
the obvious way, so that for instance a domino is north-going if and
only if it occupies a north-going domino space.  Because of symmetry,
we lose no generality by focusing on the placement probabilities
associated with north-going domino spaces.  The midpoint of the bottom
edge of each north-going domino space is some point $(\ell,m)$ with
$|\ell|+|m| \leq n-1$.  We call this the {\sl location} of the
north-going domino space.  Normalizing by dividing by $n$, we obtain
some point $(x,y)$ with $|x|+|y| < 1$.  We call this the {\sl
normalized location} of the north-going domino space.

\begin{thm}
\label{main}
Let $U$ be an open set containing the points $(\pm\frac12, \frac12)$.
If $(x,y)$ is the normalized location of a north-going domino space in
the Aztec diamond of order $n$, and $(x,y) \not\in U$, then, as $n
\rightarrow \infty$, the placement probability at $(x,y)$ is within
$o(1)$ of ${\mathcal P}(x,y)$, where
$$
{\mathcal P}(x,y) = 
\begin{cases}
0&\hbox{if $x^2+y^2\ge\frac12$ and $y<\frac12$,}\\
1&\hbox{if $x^2+y^2\ge\frac12$ and $y>\frac12$, and}\\
\frac{1}{2}+\frac{1}{\pi}\tan^{-1}
\left(\frac{2y-1}{\sqrt{1-2x^2-2y^2}}\right)
&\hbox{if $x^2+y^2<\frac12$}.\\
\end{cases}
$$
The $o(1)$ error bound is uniform in $(x,y)$ (for $(x,y) \not\in U$).
\end{thm}

Similarly, the south-going, east-going, and west-going placement
probabilities near $(x,y)$ are approximated by ${\mathcal P}(-x,-y)$,
${\mathcal P}(-y,x)$, and ${\mathcal P}(y,-x)$, respectively.  This
follows from \thmref{main} by rotational symmetry.

The organization of the rest of this article is as follows.  

In the remainder of Section~\ref{sec-introduction}, we discuss some
qualitative features of the main theorem and give some preliminaries
for the proof.  In Section~\ref{sec-creation}, we use the saddle point
method to derive asymptotic estimates for certain numbers known as
creation rates, which give placement probabilities when summed
appropriately.  In Section~\ref{sec-placement}, we use this result to
estimate, modulo an error term, the north-going placement
probabilities.  In Section~\ref{sec-exponential}, we use techniques
{}from the theory of exponential sums to justify our bound for the error
term.  This completes the proof of the theorem away from the boundary
of the diamond; Section~\ref{sec-conclusion}~provides the final
arguments that handle locations near the boundary.

Section~\ref{sec-consequences}~discusses some consequences of the
theorem.  In particular, by taking a detour through the theory of
domino tilings in general, we show that some consequences of the
arctangent formula apply not only to the particular shape we call the
Aztec diamond but also to slightly deformed versions of this shape
(\propref{robustness}).  We also give a new proof of the arctic circle
theorem and a large deviation estimate for certain properties of
random tilings of simply-connected finite regions
(\thmref{variancethm} and \propref{large}). 
Section~\ref{sec-further}~briefly sketches how the method of proof of
Theorem~\ref{main} can be adapted to handle the more general case of
random domino tilings when there is a bias in favor of one domino
orientation over the other (horizontal versus vertical).  We conclude
in Section~\ref{sec-speculations}~with speculations and open
questions.

For a treatment of the probabilistic preliminaries needed for
Section~\ref{sec-consequences}, see \cite{durrett}.

\subsection{Features of the result.}
\label{ssec-features}

As a first comment on the qualitative features of this formula, we
point out the continuity of the formula for ${\mathcal P}(x,y)$
(except at $(\pm\frac12,\frac12)$).  Indeed, if we had been so naive
as to ask for an asymptotic formula for the placement probabilities
for {\sl all} horizontal domino spaces in an asymptotically small
patch of the Aztec diamond (south-going as well as north-going), we
would not get a single value at all but rather a pair of values,
namely ${\mathcal P}(x,y)$ and ${\mathcal P}(-x,-y)$, which are not in
general equal.  That is, the local statistics are not even
approximately invariant under translations that exchange the two
color-classes.  It is therefore all the more pleasant that the local
statistics are asymptotically invariant under translations that
preserve the two color-classes (at least, they are invariant if, in
discussing local statistics, we confine ourselves to placement
probabilities, and do not inquire about correlations between
placements).

Another important feature of the formula is the singular behavior that
occurs near the normalized locations $(\pm \frac12, \frac12)$, which
we can explain as follows.  In \cite{alternating} it is shown that the
Aztec diamond of order $n$ has exactly $2^{n(n+1)/2}$ domino tilings,
and a formula derived in that article (formula~(7) of Section~4) can
be used to show that for $0 \leq k \leq n$, exactly $\binom{n}{k}
2^{n(n-1)/2}$ of the tilings have horizontal dominos covering the
leftmost squares in the first $k$ rows from the top and have vertical
dominos covering the leftmost squares in the next $n-k$ rows.  Thus,
the placement probability associated with the leftmost north-going
domino space in the $k$th row is exactly the sum
\[ 2^{-n} \sum_{i=k}^n \binom{n}{i}.\]
This sum is very close to 1 for $k -\frac{n}{2} \ll - \sqrt{n}$ and
very close to 0 for $k -\frac{n}{2} \gg \sqrt{n}$; macroscopically
speaking, the placement probability jumps from 1 to 0 discontinuously.
It might be possible to analyze the limiting behavior of the placement
probabilities in the vicinity of the singularities under suitable
scaling, but we do not explore this possibility here.

An easily-understood symmetry property of ${\mathcal P}(\cdot,\cdot)$
is the fact that
\begin{equation} \label{reflection}
{\mathcal P}(x,y) = {\mathcal P}(-x,y).  
\end{equation}
This is a consequence of the fact that reflecting a domino tiling
through the line $x=0$ carries north-going domino spaces to
north-going domino spaces.  A further identity satisfied by ${\mathcal
P}(\cdot,\cdot)$ is the relation
\begin{equation} 
\label{rotation}
{\mathcal P}(x,y) + {\mathcal P}(-y,x) + {\mathcal P}(-x,-y) +
{\mathcal P}(y,-x) = 1.
\end{equation}
To see why this is true, one need only observe that the four domino
spaces that contain a particular lattice square (fewer, if the square
is on the boundary) must have placement probabilities that sum to 1.

\begin{figure}
\begin{center}
\PSbox{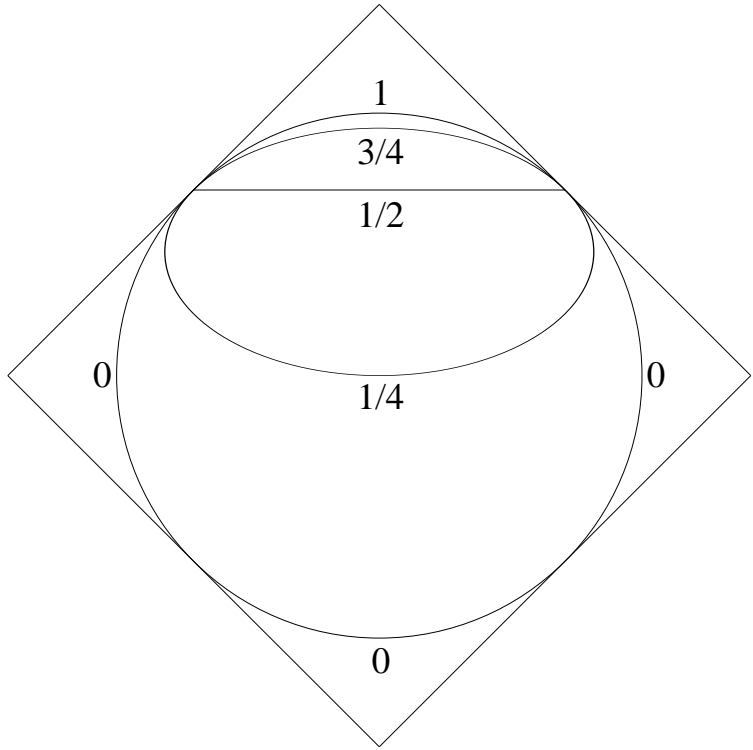 hoffset=-25 voffset=-20}{3.8in}{3.8in}
\end{center}
\caption{Level curves of north-going placement probabilities.}
\label{fig-ellipses}
\end{figure}

A subtler consequence of Theorem~\ref{main} is the fact that the level
sets of ${\mathcal P}(x,y)$ (for probabilities strictly between 0 and
1) are arcs of ellipses.  More specifically, for $0 < p < 1$ the level
set $\{(x,y): {\mathcal P}(x,y)=p\}$ and the level set $\{(x,y):
{\mathcal P}(x,y) =1-p\}$, together with the singular points $(\pm
\frac12, \frac12)$, form an ellipse tangent to the boundary of the
diamond at the two singular points.  As $p \rightarrow 0$ (or $p
\rightarrow 1$), the ellipse becomes the inscribed circle, which is
the zero-set of the function $2x^2+2y^2-1$; in the case $p=\frac12$,
the ellipse degenerates into the line segment joining the two singular
points, which is the part of the zero-set of the function $(2y-1)^2$
lying inside the Aztec diamond; and in general, the ellipse will be
the zero-set of some convex combination of $2x^2+2y^2-1$ and
$(2y-1)^2$.  The point $(0,0)$ lies on the level set $p=\frac14$,
which is an arc of an ellipse; the complementary arc of the ellipse is
the level set $p=\frac34$, and the point on this arc opposite $(0,0)$
is the point $(0,\frac23)$.  The situation is depicted schematically
in Figure~\ref{fig-ellipses}.

\subsection{Preparation for the proof.}
\label{ssec-preparation}

Recall that, under the original (unnormalized) coordinate system, each
north-going domino space in an Aztec diamond of order $n$ is assigned
some location $(\ell,m)$ with $|\ell|+|m| \leq n-1$.  It is easy to
check that $\ell+m$ must have the same parity as $n-1$.  Define
$\Pl(\ell,m;n)$ as the probability that a random domino tiling of the
Aztec diamond of order $n$ will have a domino occupying the
north-going domino space at location $(\ell,m)$; for $|\ell|+|m| >
n-1$, or $\ell+m \not \equiv n-1 \pmod{2}$, define $\Pl(\ell,m;n)=0$.
For instance, we have $\Pl(0,0;1)=\frac12$, $\Pl(0,1;2)=\frac34$, and
$\Pl(0,-1;2)=\Pl(1,0;2)=\Pl(-1,0;2)=\frac14$.

Define 
\begin{equation}
\label{creation}
\Cr(\ell,m;n) = 2(\Pl(\ell,m;n) - \Pl(\ell,m-1;n-1)).
\end{equation}
This quantity is called the {\sl net creation rate} at location
$(\ell,m)$, but the reason for this name and the interpretation of the
quantity in terms of domino shuffling are not needed for our purposes.
(For the motivation, see \cite{gip}.)

Define $c(a,b;n)$ to be the coefficient of $z^a$ in
$(1+z)^{n-b}(1-z)^b$.  (Note that $c(a,b;n)$ is the Krawtchouk
polynomial $P_a$ evaluated at $b$.  For information about Krawtchouk
polynomials, see \cite[p.~130]{codes}.)  Our proof of \thmref{main}
will be based on the following result from \cite{gip}:

\begin{prop}
\label{gessel}
Let $n>0$.  Suppose $\ell$ and $m$ are integers with $\ell+m \equiv n
\pmod{2}$ and $|\ell|+|m| \le n$.  If we let $a=(\ell+m+n)/2$ and
$b=(\ell-m+n)/2$, then
$$
\Cr(\ell,m;n+1)=c(a,b;n)c(b,a;n)/2^n.
$$
For other integers $\ell$ and $m$, we have $\Cr(\ell,m;n+1)=0$.
\end{prop}

This proposition implies that the creation rates are non-negative, if
we use the identity $c(b,a;n)b!(n-b)!=c(a,b;n)a!(n-a)!$.  (This
identity is a standard fact about Krawtchouk polynomials, and follows
immediately from Theorem~17 on page~152 of \cite{codes}.)  When
combined with \propref{gessel}, the identity implies that
$\Cr(\ell,m;n+1)$ is a positive factor times the square of a
Krawtchouk polynomial, and hence that $\Cr(\ell,m;n+1) \ge 0$.  Note
that this inequality, together with (\ref{creation}), yields
\begin{equation}
\label{incbound}
\Pl(\ell,m;n) \le \Pl(\ell,m+h;n+h)
\end{equation}
for $h \ge 0$ by induction on $h$.

We will also need the following result on exponential sums.  

\begin{thm}[Kusmin-Landau]
\label{kusmin-landau}
Let $||\cdot||$ denote the distance to the nearest integer, $I$ be an
interval, and $f$ be a real-valued function on $I$.  If $f$ is
continuously differentiable, $f'$ is monotonic, and $ ||f'|| \ge
\lambda > 0$ on $I$, then
$$\sum_{n \in I\cap\Z}\exp(2\pi i \,f(n)) = O(\lambda^{-1}).$$
The constant implicit in the
$O(\lambda^{-1})$ term does not depend on $I$.
\end{thm}

\noindent A proof can be found in \cite[p.~7]{expsums}.

\subsection{Outline of the Proof of Theorem~1}
\label{ssec-outline}

We begin the proof of \thmref{main} by using (\ref{creation}) to write
the placement probabilities as sums of creation rates, which gives the
formula
\begin{equation}
\label{crsum}
\Pl(\ell,m;n) = \frac{1}{2}\sum_{k \ge 0}{\Cr(\ell,m-k;n-k)}.
\end{equation}
(Note that the remark after \propref{gessel} shows that, as claimed in
the abstract, this is a weighted sum of squares of Krawtchouk
polynomials.)  We will estimate the creation rates, and then use our
estimate to prove the asymptotic formula for the placement
probabilities.

Because of \propref{gessel}, to estimate the creation rates it
suffices to approximate the coefficients of the polynomials
$(1+z)^{n-b}(1-z)^b$.  To do this, we write the coefficients as
contour integrals in the usual way, and then apply the saddle point
method to these integrals.  Sufficiently far outside the arctic
circle, this method shows that the creation rates are exponentially
small in $n$ (\propref{expsmall}); sufficiently far inside, it
approximates them by a well-behaved function times an oscillating
factor (\propref{crest}).  With additional work, it might be possible
to obtain a uniform estimate over the entire Aztec diamond, but we can
make do with just these estimates.

We would then like to substitute our creation rate estimates into
(\ref{crsum}) and convert the sum to an integral to determine its
asymptotics.  If we are willing to be unrigorous, we can wishfully
replace the oscillatory cosine-squared factor in \propref{crest} by
its mean value $\frac12$, obtaining (for locations inside the
inscribed circle)
\begin{eqnarray*}
\Pl(\ell, m; n)  
& \approx & \frac12 \sum_{k = 0}^{t_{{\rm max}}} 
	\frac{2}{\pi \sqrt{(n-k)^2-2\ell^2-2(m-k)^2}} \\
& \approx & \frac12 \int_0^{t_{{\rm max}}}
	\frac{2}{\pi \sqrt{(n-k)^2-2\ell^2-2(m-k)^2}} \ dk \\
& = & \frac{1}{\pi} \left. \left( \tan^{-1} 
	\frac{k+n-2m}{\sqrt{(n-k)^2-2\ell^2-2(m-k)^2}} \right)
	\right|_{k=0}^{k=t_{{\rm max}}} \\
& = & \frac12 - \frac{1}{\pi} \tan^{-1}
\frac{n-2m}{\sqrt{n^2-2\ell^2-m^2}} \\ 
& = & \frac12 + \frac{1}{\pi} \tan^{-1}
\frac{2y-1}{\sqrt{1-2x^2-2y^2}}, 
\end{eqnarray*}
where $k=t_{{\rm max}}$ is the larger of the two roots of the equation
$(n-k)^2-2\ell^2-2(m-k)^2=0$, $x=\ell/n$, and $y=m/n$; we truncate the
sum and integral at $t_{{\rm max}}$ on the supposition (to be
discussed in the next paragraph) that essentially no creation occurs
outside the arctic circle.  To make this argument rigorous, we need to
deal honestly with the oscillating factor in the creation rate
estimate inside the arctic circle, and we need to circumvent the
non-uniformity of our estimates.

The non-uniformity of the estimates can be dealt with simply by
summing over a smaller interval than in (\ref{crsum}).  The
exponentially small bounds on the creation rates outside the arctic
circle show that the terms in (\ref{crsum}) that come from locations
outside the arctic circle contribute very little to the sum.
Motivated by this, we look at the sum of all the terms that come from
locations that are far enough inside the arctic circle that our
creation rate estimates from \propref{crest} apply.  Because the
creation rates are all non-negative, this new sum underestimates the
placement probability.  Dealing appropriately with the oscillating
factor (as described below) gives an estimate for the new sum; as we
see in the computation above, this estimate turns out to be the
arctangent formula from \thmref{main}.  A short argument shows that
the placement probabilities can be no larger asymptotically
(\propref{maininside}), and because our estimate is an underestimate
we know they can be no smaller.  This completes the proof.  (Actually,
this method works only away from the edges of the diamond, so it is
not until Section~\ref{sec-conclusion}~that the proof is completed.)

All that remains is to describe how to deal with the oscillating
factor in the summand.  We must show that replacing the oscillating
factor by its average value has an asymptotically negligible effect on
the sum.  Equivalently, we must show that the difference between the
original sum and the smoothed sum is small.  This difference is an
exponential sum, and we can estimate it using the Kusmin-Landau
Theorem once some preparatory results (Lemmas~\ref{phase},
\ref{algebraic}, and \ref{algroots}) are in place.

\section{Creation Rate Estimates}
\label{sec-creation}

Our proof of the asymptotic formula for placement probabilities begins
with an estimate of creation rates, which is proved using the saddle
point method.  Because \propref{gessel} is most conveniently stated
for an Aztec diamond of order $n+1$, we will estimate the creation
rates in an Aztec diamond of order $n+1$.  From this point on, we
assume that $\ell+m \equiv n \pmod{2}$, because otherwise
$\Cr(\ell,m;n+1)$ is necessarily $0$.  As pointed out in
subsection~\ref{ssec-outline}, the creation rates behave differently
inside and outside the inscribed circle.  If we estimate the creation
rates inside it, we get the following result:

\begin{prop}
\label{crest}
Fix $\varepsilon > 0$.  If $\ell^2 + m^2 \le (1-\varepsilon)n^2/2$ and
$\ell+m \equiv n \pmod{2}$, then
$$
\Cr(\ell,m;n+1) = \frac{4\cos^2
\Phi(\ell,m;n)}{\pi\sqrt{n^2-2\ell^2-2m^2}} + O_\varepsilon(n^{-2})
$$
for some function $\Phi(\ell,m;n)$, which we determine explicitly
below. 
\end{prop}

The subscript in $O_\varepsilon(n^{-2})$ indicates that the implicit
constant depends on $\varepsilon$.  In this paper, if any subscripts
appear on a big $O$ term, then the implicit constant depends only on
the indicated variables, but the absence of subscripts should not be
taken to imply that the implicit constant is absolute.

\begin{proof}
Let 
\begin{equation}
\label{defoff}
f(z) = \frac{(1+z)^{n-b}(1-z)^b}{z^a},
\end{equation}
where $a$ and $b$ are defined as in the statement of \propref{gessel}.
To approximate the creation rate, we need to approximate $c(a,b;n)$,
which is the constant term of $f(z)$.  The constant term is given by
the usual contour integral, which we will approximate using the saddle
point method.

Write $a=(1+u)n/2$ and $b=(1+v)n/2$, so that $-1 \le u,v \le 1$.  Note
that the $u,v$ coordinates are related to the coordinates in the
statement of the proposition by $u=(\ell+m)/n=x+y$ and
$v=(\ell-m)/n=x-y$.  We will keep $u$ and $v$ fixed as we send $n$ to
infinity.

The critical points of $f(z)$ are
$$z_1 = \frac{-v+\sqrt{u^2+v^2-1}}{1-u}$$
and
$$z_2 = \frac{-v-\sqrt{u^2+v^2-1}}{1-u}.$$
Because $\ell^2+m^2<n^2/2$, we have $u^2+v^2<1$.  It follows that
$z_1$ and $z_2$ are complex conjugates on the circle $|z|^2=
(1+u)/(1-u)$.  We now apply the saddle point method.  To find the
constant term of $f(z)$, we integrate $f(z)/(2\pi i \, z)$ about the
circle of radius $\sqrt{(1+u)/(1-u)}$ centered at the origin.  One can
check that, on this circle, $|f(z)|$ is greatest at the critical
points $z_1$ and $z_2$.  (To check it, parametrize the circle by the
angle $\theta$ formed with the real axis.  One has $\partial \log
|f(z)|^2/\partial \theta = 0$ iff $z$ is one of the two critical
points or $z$ lies on the real axis.  At the critical points,
$\partial^2 \log |f(z)|^2 / \partial \theta^2 = n(u^2+v^2-1)/(1-v^2) <
0$, so $|f(z)|$ has maxima at these points.  It must have minima on
the real axis, since between any two local maxima there must be a
local minimum.)  As $n$ goes to infinity, the integral is given
asymptotically by the integrals over the parts of the path near the
critical points, which can be estimated straightforwardly.  This is
the saddle point method.  We will omit the details of the argument
leading to the approximation, because they are standard, and can be
found, for example, in \cite[pp.~87--89]{asymp}.

The saddle point method tells us that the constant term of $f(z)$ is
the sum $Z_1(1+O(n^{-1}))+Z_2(1+O(n^{-1}))$, where
\begin{equation}
\label{defofz1}
Z_1 = \frac{f(z_1)}{2\pi  z_1}
\sqrt{\frac{2\pi}{(\log f)''(z_1)}}
\end{equation}
and 
\begin{equation}
\label{defofz2}
Z_2 = {\overline{Z_1}} = \frac{f(z_2)}{2\pi 
z_2}\sqrt{\frac{2\pi}{(\log f)''(z_2)}}.
\end{equation} 
(For the proof of \propref{crest} we will not need to determine the
signs of the square roots in \eqref{defofz1} and \eqref{defofz2}, but
they must be chosen so that $Z_1$ and $Z_2$ are complex conjugates.)

Simplifying $z^2(\log f)''(z)$ yields
$$
z^2(\log f)''(z) = {\frac {n\left
(1+u-4\,z^{2}-2\,uz^{2}-4\,vz^{3}+uz^{4}-z^{4}\right )}
{2\,(z^2-1)^{2}}}.
$$
{}From this, one can check that at either critical point of $f(z)$,
$z^2(\log f)''(z)$ has absolute value
\begin{equation}
\label{absofl}
|z_i^2(\log f)''(z_i)|
=
\frac{n}{2}\sqrt{\frac{(1-u^2-v^2)(1-u^2)}{(1-v^2)}}.
\end{equation}

Let $\Psi(u,v;n)$ be the phase of $Z_1$, so that $Z_1 = |Z_1|
\exp(i\,\Psi(u,v;n)).$ Then
$$Z_1+Z_2 = 2|Z_1|\cos\Psi(u,v;n),$$
and $c(a,b;n)$ is approximated by
\begin{equation}
\label{approx}
c(a,b;n)= 
2|Z_1|\cos\Psi(u,v;n)+O\left(\frac{|Z_1|}{n}\right) .
\end{equation}
Of course,
\begin{equation}
\label{absofz1}
|Z_1| = \frac{|f(z_1)|}{2\pi}\sqrt{\frac{2\pi}{|z_1^2(\log
f)''(z_1)|}}.
\end{equation}
Since $|1+z_1|^2 = 2(1-v)/(1-u)$, $|1-z_1|^2=2(1+v)/(1-u)$, and
$|z_1|^2=(1+u)/(1-u)$, we see that
\begin{equation}
\label{absoff}
|f(z_1)| =
2^{n/2}\frac{(1-v)^{(n-b)/2}(1+v)^{b/2}}{(1-u)^{(n-a)/2}(1+u)^{a/2}}.
\end{equation}

We have $\Cr(\ell,m;n+1) = c(a,b;n)c(b,a;n)/2^n,$ by \propref{gessel}.
Interchanging $a$ and $b$ corresponds to interchanging $u$ and $v$.
Let $\tilde x$ denote the result of interchanging $u$ and $v$ (and
also $a$ and $b$) in the expression $x$, so that, for example,
$\widetilde{u-2v} = v-2u$.  When we substitute (\ref{absoff}) and
(\ref{absofl}) into (\ref{absofz1}), we see that
$$
|Z_1||\widetilde{Z_1}| = \frac{2^n}{\pi n \sqrt{1-u^2-v^2}}.
$$
Hence, by (\ref{approx})
$$
\Cr(\ell,m;n+1)=
\frac{4\cos\Psi(u,v;n)\cos\Psi(v,u;n)}{\pi n
\sqrt{1-u^2-v^2}}+O(n^{-2}).
$$
(To see that the error term is $O(n^{-2})$, one uses the fact that it
is $O(|Z_1||\widetilde{Z_1}|/(n2^n))$ and that $|Z_1||\widetilde{Z_1}|
= O(2^n/n)$.)

Now we check that $\cos\Psi(v,u;n)=\pm\cos\Psi(u,v;n).$ The identity
$$
c(b,a;n)b!(n-b)!=c(a,b;n)a!(n-a)!
$$
suggests that this should be so, but does not seem to prove it.  If we
set $\alpha = z_1^2(\log f)''(z_1)$, we find (after some
simplification) that
$$
\frac{\alpha}{\tilde \alpha} = \frac{1-u^2}{1-v^2}.
$$
Thus, $\alpha$ and $\tilde \alpha$ have the same phase.  If we combine
the formulas
$$
\frac{1+z_1}{1+{\tilde z_1}} = \frac{1-v}{1-u}
$$
and
$$
\frac{1-z_1}{1+z_1} = \frac{1+v}{1-v}\,\frac{1}{\tilde z_1}
$$
with
$$
f(z_1) =
(1+z_1)^n\left(\frac{1-z_1}{1+z_1}\right)^b
\left(\frac{1}{z_1}\right)^a,
$$
we find that $f(z_1)$ equals $\widetilde{f(z_1)}$ times a positive
factor, so their phases are equal.  Because
$$
Z_1 = \pm \frac{f(z_1)}{2\pi}\sqrt{\frac{2\pi}{\alpha}}
\quad \hbox{and} \quad
\widetilde{Z_1} = \pm
\frac{\widetilde{f(z_1)}}{2\pi}\sqrt{\frac{2\pi}{\tilde \alpha}},
$$
we see that $Z_1$ and $\widetilde{Z_1}$ have the same phase, to within
a sign, so $\cos\Psi(v,u;n)=\pm\cos\Psi(u,v;n).$

Finally, we change to the coordinates of our generating function by
the substitutions $u=(\ell+m)/n$ and $v=(\ell-m)/n$.  We set
$\Phi(\ell,m;n) = \Psi(u,v;n)$.  Then when $\ell+m\equiv n \pmod{2}$,
the creation rate at the $(\ell,m)$ location in an Aztec diamond of
order $n+1$ is
$$\Cr(\ell,m;n+1) =
\frac{\pm4\cos^2\Phi(\ell,m;n)}{\pi \sqrt{n^2-2\ell^2-2m^2}}+O(n^{-2}).
$$
Because creation rates must be non-negative, the $\pm$ sign in this
formula can always be taken to be $+$.

The constant implicit in the big $O$ depends continuously on $u$ and
$v$.  Thus, for fixed $\varepsilon > 0$ the constant can be chosen
uniformly for all $u$ and $v$ with $u^2+v^2\le1-\varepsilon$.  We have
therefore proved the result claimed in the statement of the
proposition.
\end{proof}

Given $\ell$, $m$, and $n$, define
\begin{equation}
\label{sdef}
S_\varepsilon = \{k \in \Z : k \ge 0 \hbox{ and }
\ell^2+(m-k)^2\le(1-\varepsilon)(n-k)^2/2 \}.
\end{equation}
Also, define $k_{{\rm max}}$ to be the greatest element of
$S_\varepsilon$, and $k_{{\rm min}}$ to be the least.  (In
Section~\ref{sec-placement}, we will sum the creation rates in
(\ref{crsum}) as $k$ varies over $S_\varepsilon$.  We will do so to
make it possible to apply \propref{crest} to the terms in the sum, as
described in subsection~\ref{ssec-outline}.)

\begin{lem}
\label{maxbound}
Suppose that $|\ell|+|m|
\le(1-\delta)n$ for some fixed $\delta > 0$.
Then
$k_{{{\rm max}}}
\le
(1-(2-\sqrt{2})\delta + O(\varepsilon))n$.  Hence, 
if $\varepsilon$ is small enough compared to $\delta$, then
$(n-k_{{{\rm max}}})^{-1} =
O(n^{-1}),$
and for $k \in S_\varepsilon$ we have
$$
\frac{1}{\sqrt{(n-k)^2-2\ell^2-2(m-k)^2}} =
O(\varepsilon^{-1/2}n^{-1}) 
 = O_{\varepsilon,\delta}(n^{-1}).
$$
\end{lem}

\begin{proof}
We have $k_{{{\rm max}}} =
(2y-1+\sqrt{2((1-y)^2-x^2)}+O(\varepsilon))n,$ where $x=\ell/n$ and
$y=m/n$.  For fixed $y$, this function is clearly maximized at $x=0$.
When $x=0$, it becomes a linear function maximized at $y=1-\delta$
(for $|x|+|y|\le 1-\delta$).  This yields $k_{{{\rm max}}} \le
(1-(2-\sqrt{2})\delta+O(\varepsilon))n$.  Therefore, $(n-k_{{{\rm
max}}})^{-1}=O(n^{-1})$ if $\varepsilon$ is small enough compared to
$\delta$.

For $k \in S_\varepsilon$, we have
$$
\ell^2+(m-k)^2 \le (n-k)^2(1-\varepsilon)/2
.$$
It follows that
$$
(n-k)^2-2\ell^2-2(m-k)^2 \ge (n-k)^2\varepsilon.
$$
Therefore,
$$
\frac{1}{\sqrt{(n-k)^2-2\ell^2-2(m-k)^2}} =
O(\varepsilon^{-1/2}(n-k)^{-1})  = O(\varepsilon^{-1/2}n^{-1}).
$$
If we are not worrying about dependence on $\varepsilon$, this is
$O(n^{-1})$.
This completes the proof.
\end{proof}

In order to apply the Kusmin-Landau Theorem to deal with the
exponential sums that will appear later in the proof (as discussed in
subsection~\ref{ssec-outline}), we will need to specify $\Phi$, since
the phase of $Z_1$ is not uniquely determined.  Also, it will be
convenient to extend it to a function of real, and even complex,
variables (rather than just integers).

Given a point $(x,y) \ne (0,0)$ in the plane, define $\theta(x,y)$ to
be the angle in $(-\pi,\pi]$ formed by the right half of the
horizontal axis and the ray from the origin through $(x,y)$.

\begin{lem}
\label{phase}
We can choose $\Phi(\ell,m;n)$ in \propref{crest} so that 
if one sets
$\ell=xn$, $m=yn$, and $k={\kappa}n$ in
$d\Phi(\ell,m-k;n-k)/dk$, then $d\Phi(\ell,m-k;n-k)/dk$ equals
\begin{eqnarray*}
\theta\left({-x+y-{\kappa}},
{\sqrt {(1-{\kappa})^2-2\,(x^2+(y-{\kappa})^2)}
}\right) - \\
\theta\left({1-{\kappa}-2\,x},
{\sqrt {(1-{\kappa})^2-2\,(x^2+(y-{\kappa})^2)}
}\right) + \\
{\frac {{x}^{2}-{\kappa}-3\,y{\kappa}+2\,{{\kappa}}^{2}+y+{y}^{2}}
{n\sqrt{(1-{\kappa})^2-2\, (x^2+(y-{\kappa})^2)}\,
(y+1-2\,{\kappa}-x)(y+1-2\,{\kappa}+x)}}.
\end{eqnarray*}
As $n$ tends to infinity, the last term is $O(1/n)$ for $k \in
S_\varepsilon$ with $\varepsilon > 0$ fixed, as long as $|\ell|+|m|
\le (1-\delta)n$ for some fixed $\delta > 0$, and $\varepsilon$ is
small enough compared to $\delta$.
\end{lem}

\begin{proof}
{}From (\ref{defofz1}), we see that we can
choose $\Phi(\ell,m;n)$ to be the imaginary part 
\begin{equation}
\label{choice}
\hbox{Im}
\left (
\log f(z_1)  - \log z_1 - \frac{1}{2}\log((\log f)''(z_1))
\right)
.
\end{equation}

If we substitute $m-k$ for $m$ and $n-k$ for $n$ and differentiate,
then the first term of (\ref{choice}) contributes the $\theta$-terms
in the formula we are proving.  To see this, recall that (up to an
irrelevant multiple of $2\pi i$)
$$
\log f(z_1) = (n-b)\log(1+z_1) + b \log(1-z_1) - a\log z_1.
$$
After we express this in terms of $n$, $\ell$, and $m$ and substitute
$m-k$ for $m$ and $n-k$ for $n$, the right hand side becomes
\begin{equation}
\label{messyl}
(n/2-k-\ell/2+m/2)\log(1+{\hat z_1}) + 
(n/2+\ell/2-m/2) \log(1-{\hat z_1}) \cr
- (n/2-k+\ell/2+m/2)\log {\hat z_1 },
\end{equation}
where $\hat z_1$ is the function of $k$ that results from making the
substitutions in $z_1$.  Denote by $L$ the function (\ref{messyl}).
When we differentiate $L$ with respect to $k$ (holding $n$, $\ell$,
and $m$ fixed), we get
$$
\frac{\partial L}{\partial k} = 
\log{\hat z_1} - \log(1+{\hat z_1}) + 
\frac{\partial L}{\partial {\hat z_1}}
\frac{\partial {\hat z_1}}{\partial k}.
$$
Because $z_1$ is a critical point of $f$, ${\partial L}/{\partial
{\hat z_1}} = 0$, so ${\partial L}/{\partial k} = \log{\hat z_1} -
\log(1+{\hat z_1})$.  Now expressing the imaginary parts of the
logarithms in terms of $\theta$ gives the desired terms from the
formula we are proving.  (To simplify the terms to the form found in
the statement of the lemma, one has to use the fact that for $\alpha >
0$, $\theta(\alpha x,\alpha y) = \theta(x,y)$.)

When we substitute and differentiate, the remaining terms in
(\ref{choice}) clearly give algebraic results.  We omit the details of
the calculations, since they are tedious and straightforward.

The claim that the last term is $O(1/n)$ for $k \in S_\varepsilon$ is
a consequence of \lemref{maxbound}.  The only thing to check is that
although the denominator vanishes at ${\kappa} = (y+1\pm x)/2$, these
two points are never in $S_\varepsilon$ (or near enough to cause
problems).  To see that, note that the definition (\ref{sdef}) of
$S_\varepsilon$ is equivalent to the set of $k \ge 0$ for which
\begin{equation}
\label{othersdef}
(1-\kappa)^2-2(x^2+(y-\kappa)^2) \ge
\varepsilon(1-\kappa)^2.
\end{equation}
Note that $\kappa = 1$ is impossible (since then we must have $x=0$
and $y=1$, so $|x|+|y| > 1-\delta$).  However, substituting $\kappa =
(y+1\pm x)/2$ in the left hand side of (\ref{othersdef}) gives $-(3x
\pm(1-y))^2/4$.  Thus, the factors $y+1-2\kappa -x$ and
$y+1-2\kappa+x$ in the denominator of the last term in our main
formula cannot become arbitrarily small, and the last term is indeed
$O(1/n)$.
\end{proof}

In Section~\ref{sec-exponential}, we will need the following result
(to make it possible to apply exponential sum techniques in the way
described in subsection~\ref{ssec-outline}).

\begin{lem}
\label{algebraic}
The function $d^2 \Phi(\ell,m-k;n-k)/dk^2$ is algebraic.  For any
fixed $n$, $\ell$, and $m$ satisfying $|\ell|+|m| < n$ and
$\varepsilon > 0$, there exists a neighborhood $U$ in $\C$ of the
smallest real interval containing $S_\varepsilon$ such that as a
function of $k$, $d^2 \Phi(\ell,m-k;n-k)/dk^2$ is holomorphic on $U$.
\end{lem}

\begin{proof}
We will use the formula for $d \Phi(\ell,m-k;n-k)/dk$ from
\lemref{phase}.  Let $U$ be a small, simply-connected neighborhood in
$\C$ of the smallest real interval containing $S_\varepsilon$, such
that the points $k=n(y+1\pm x)/2$ are not in $U$.  (We checked at the
end of the proof of \lemref{phase} that these points are not in
$S_\varepsilon$.)  It follows from the definition of $S_\varepsilon$
that
$$
(n-k)^2-2(\ell^2+(m-k)^2) \ge \varepsilon(n-k)^2 \ge 0
$$
on $S_\varepsilon$.  If $n=k$, then $\ell = 0$ and $m=k=n$, so
$|\ell|+|m| = n$ (contradicting $|\ell|+|m|<n$).  Thus,
$(n-k)^2-2(\ell^2+(m-k)^2) \ge \varepsilon$ on $S_\varepsilon$, and
hence there is a holomorphic square root of
$(n-k)^2-2(\ell^2+(m-k)^2)$ on $U$, if $U$ was chosen to be
sufficiently small.  It follows that the third term (the algebraic
term) of the formula for $d \Phi(\ell,m-k;n-k)/dk$ in \lemref{phase}
is holomorphic on $U$.  The derivative of that term is thus algebraic
and holomorphic on $U$, so to complete the proof we just need to check
this for the other two terms.

The first two terms can be expressed in terms of the arctangent.  If
we do so, we find that the derivative with respect to $k$ of the sum
of those two terms is
$$
\frac{-3x^2+2y\kappa+1-2\kappa-y^2}
{n\sqrt{(1-\kappa)^2-2(x^2+(y-\kappa)^2)}
(y+1-2\kappa-x)(y+1-2\kappa+x)}.
$$ 
This is also algebraic and holomorphic on $U$.  Thus,
$d^2\Phi(\ell,m-k;n-k)/dk^2$ is holomorphic on $U$ and algebraic, as
desired.
\end{proof}

We now know everything we need to know about how the creation rates
behave inside the arctic circle.  Outside the arctic circle, we can
get an exponentially small upper bound for the creation rates.  This
will be used for bounding the placement probabilities outside the
arctic circle (\propref{arcticbound}).

\begin{prop}
\label{expsmall}
For each $\varepsilon > 0$, there exists a positive constant $r < 1$
such that whenever $\ell^2+m^2 > (1+\varepsilon)n^2/2$,
$$
\Cr(\ell,m;n+1) = O(r^n).
$$
\end{prop}

\begin{proof}
We assume that $\ell+m\equiv n \pmod{2}$, since otherwise
$\Cr(\ell,m;n+1)=0$.  As in the proof of \propref{crest}, we will
integrate $f(z)/(2\pi i \, z)$ around a circle about the origin,
where, as in (\ref{defoff}),
$$
f(z) = \frac{(1+z)^{n-b}(1-z)^b}{z^a}.
$$
However, since we are looking only for an upper bound and not for an
asymptotic estimate, we will not need the full saddle point method.
We will only sketch the proof, because the details are straightforward
but somewhat tedious to check.

We will use the same notation as in the proof of \propref{crest}; for
example, we write $a=(1+u)n/2$ and $b=(1+v)n/2$.  Since
$u^2+v^2>1+\varepsilon$, the critical points
$$
z_1 = \frac{-v+\sqrt{u^2+v^2-1}}{1-u}
$$ 
and
$$
z_2 = \frac{-v-\sqrt{u^2+v^2-1}}{1-u}
$$ 
of $f(z)$ are real.  (Of course, the case $u=1$ has to be handled
separately, but this will not cause problems.)  We will integrate
$f(z)/(2\pi i \, z)$ around a circle of radius $R$, where $R$ will be
either $|z_1|$ or $|z_2|$.  We choose $R = |z_i|$ where $|f(z_i)|$ is
the lesser of $|f(z_1)|$ and $|f(z_2)|$.  To bound the integral, we
will use the fact that the absolute value of the integral is at most
as large as the greatest value of $|f(z)|$ on the circle.

It is not hard to check by straightforward manipulation of
inequalities that $|f(z_1)| > |f(z_2)|$ if $uv > 0$, and $|f(z_1)| <
|f(z_2)|$ if $uv < 0$.  (Since $|u|,|v| \le 1$ and $u^2+v^2 >
1+\varepsilon$, we cannot have $uv=0$.)  Thus, $R=|z_2|$ if $uv >0$,
and $R=|z_1|$ otherwise.

Take $i \in \{1,2\}$ so that $R=|z_i|$.  On the circle of radius $R$
about $0$, $|f(z)|$ is greatest when $z=z_i$; in fact, the second
derivative test shows that this is the only local maximum.  Thus, the
integral is bounded by $|f(z_i)|$, so $|c(a,b;n)| \le |f(z_i)|$.

Because the sign of $uv$ doesn't change when $u$ and $v$ are
interchanged, $|\widetilde{f(z_i)}|$ is the lesser of
$|\widetilde{f(z_1)}|$ and $|\widetilde{f(z_2)}|$.  Hence, $|c(b,a;n)|
\le |\widetilde{f(z_i)}|$.  It follows that
$$
\Cr(\ell,m;n+1) \le \frac{|f(z_i)||\widetilde{f(z_i)}|}{2^n}.
$$

A simple calculation gives $|f(z_1)||\widetilde{f(z_2)}| =
|f(z_2)||\widetilde{f(z_1)}| = 2^n$.  The inequalities $|f(z_i)| <
|f(z_{3-i})|$ and $|\widetilde{f(z_i)}| < |\widetilde{f(z_{3-i})}|$,
together with the fact that the only dependence on $n$ in any of these
expressions is in the exponent, imply that the creation rate at
$(u,v)$ is $O(r^n)$ for some $r < 1$.  A little more care in the
estimates shows that this bound can be chosen uniformly for $u^2+v^2 >
1+\varepsilon$, as desired.
\end{proof}

\section{Placement Probability Estimates}
\label{sec-placement}

Now that we know the creation rates, we can determine the placement
probabilities.  Fix $\delta > 0$.  In this section, we will look at
the placement probabilities $\Pl(\ell,m;n+1)$ at points $(\ell,m)$
satisfying $\ell+m \equiv n \pmod{2}$ and $|\ell| + |m| \le
(1-\delta)n$.  (The congruence condition rules out the placement
probabilities that we know are 0, and the inequality lets us apply
results such as Lemmas~5--7.)

{}From (\ref{crsum}), we see that
$$
\Pl(\ell,m;n+1) = \frac{1}{2}\sum_{k \ge 0}{\Cr(\ell,m-k;n+1-k)}.
$$
It will turn out that the creation rates on or beyond the arctic
circle make a vanishing contribution to this sum as $n \rightarrow
\infty$, so we can remove them from the sum without affecting its
asymptotics.  To remove these terms from the sum, fix $\varepsilon >
0$ (which we assume is small compared to $\delta$, so that we can
apply results such as \lemref{phase}), and look at the sum
$$
\Plt_\varepsilon = \frac{1}{2}\sum_{k \in S_\varepsilon}{\Cr(\ell,m-k;n+1-k)},
$$
where $S_\varepsilon$ is defined by (\ref{sdef}).  (Note that
sometimes $S_\varepsilon = \emptyset$; we will see that this occurs
only when the placement probability is nearly $0$.)  We will
approximate $\Plt_\varepsilon$, and prove that it approximates
$\Pl(\ell,m;n+1)$.  First, we prove a few easy lemmas.

\begin{lem}
\label{bound}
Consider the equation $(1-t)^2-2x^2-2(y-t)^2=0$.  For $|x|+|y|<1$,
this equation has two real roots $t$.  The greater root is $0$ iff
$x^2+y^2=1/2$ and $y < 1/2$, and is less than $0$ iff $x^2+y^2>1/2$
and $y < 1/2$.  The lesser root is greater than or equal to $0$ iff
$x^2+y^2 \ge 1/2$ and $y > 1/2$.
\end{lem}

\begin{proof}
Since the discriminant of the polynomial is $8(1-x-y)(1+x-y)$, we see
that it has two real roots whenever $|x|+|y|<1$.  Clearly, $0$ is a
root iff $x^2+y^2=1/2$, and since the sum of the roots is $4y-2$, it
is the greater root iff also $y < 1/2$.  One can check the other
claims similarly.
\end{proof}

\begin{lem}
\label{estim}
Let $\delta > 0$, and suppose $|x|+|y| \le 1-\delta$.  Let
$\kappa(\varepsilon)$ be any branch of the multivalued function of
$\varepsilon$ defined by $(1-\varepsilon)(1-\kappa)^2 - 2x^2 -
2(y-\kappa)^2 = 0$.  Then for $\varepsilon > 0$ (and sufficiently
small relative to $\delta$), we have $\kappa(\varepsilon) = \kappa(0)
+ O_\delta(\varepsilon)$.  (The constant implicit in the
$O_\delta(\varepsilon)$ does not depend on $x$, $y$, or
$\varepsilon$.)
\end{lem}

\begin{proof}
This simply amounts to showing that $\kappa'(\varepsilon)$ is bounded
as a function of $x$, $y$, and $\varepsilon$, for $\varepsilon$
sufficiently small.  If one computes $\kappa(\varepsilon)$ using the
quadratic formula, and then differentiates it with respect to
$\varepsilon$, one finds that it equals a continuous function of $x$,
$y$, and $\varepsilon$ (for $\varepsilon$ near $0$) divided by
$$
\sqrt{(1-\varepsilon)(y-1)^2-(1+\varepsilon)x^2}.
$$
If $\varepsilon$ is small enough compared to $\delta$, then
$\kappa'(\varepsilon)$ will be continuous, and hence bounded, for all
$x$, $y$, and $\varepsilon$ with $|x|+|y|\le 1-\delta$.  Then
$\kappa(\varepsilon)=\kappa(0)+O(\varepsilon)$, as desired.
\end{proof}

\begin{prop}
\label{firstapprox}
Let $\delta > 0$ and $\varepsilon > 0$, such that $\varepsilon$ is
sufficiently small compared to $\delta$.  If $|\ell| + |m| \le
(1-\delta)n$ and $\ell+m \equiv n \pmod{2}$, then
$$
{\Plt_\varepsilon} = {\mathcal P}(\ell/n,m/n)
+ O_{\varepsilon, \delta}(n^{-1}) + O_\delta(\varepsilon^{1/2}).
$$
\end{prop}

\begin{proof}
{}From formula (\ref{crsum}) and \propref{crest}, we see that
$\Plt_\varepsilon$ is approximated by $$\frac{1}{2}\sum_{k \in
S_\varepsilon}{\frac{4\cos^2\Phi(\ell,m-k;n-k)}{\pi
\sqrt{(n-k)^2-2\ell^2-2(m-k)^2}}}+O(n^{-1}).$$ (The error here is
bounded by a constant multiple of $\sum_{j\ge n-k_{{\rm max}}}j^{-2}$,
which is $O((n-k_{{\rm max}})^{-1}),$ and hence $O(n^{-1})$ by
\lemref{maxbound}.)

Since $4\cos^2 z = \exp(2iz) + \exp(-2iz) + 2$, we see that
$\Plt_\varepsilon$ is given to within $O(n^{-1})$ by the sum of
\begin{equation}
\label{firstsum}
\frac{1}{2}\sum_{k \in S_\varepsilon}{\frac{2}{\pi
\sqrt{(n-k)^2-2\ell^2-2(m-k)^2}}},
\end{equation}
with
$$
\frac{1}{2}\sum_{k \in S_\varepsilon}{\frac{1}{\pi
\sqrt{(n-k)^2-2\ell^2-2(m-k)^2}}} \exp(2i\,\Phi(\ell,m-k;n-k))
$$
and its complex conjugate.  \propref{expvanish} will show that the
latter two sums are $O(n^{-1})$ as $n$ goes to infinity.

Assuming \propref{expvanish}, we can prove the desired limit by
approximating the sum (\ref{firstsum}) with an integral.  The sum is
equal to
\begin{equation}
\label{intap}
\frac{1}{2}\int_{k_{{\rm min}}}^{k_{{\rm max}}} {\frac{2}{\pi
\sqrt{(n-k)^2-2\ell^2-2(m-k)^2}}}\,dk 
+O(n^{-1});
\end{equation}
to see that the error is $O(n^{-1})$, note that the summand (viewed as
a function of a real variable $k$) is $O(n^{-1})$ by
\lemref{maxbound}, and is monotonic on $S_\varepsilon \cap
(-\infty,2m-n)$ and $S_\varepsilon \cap (2m-n,\infty)$.

By \lemref{bound}, the polynomial $(n-t)^2-2\ell^2-2(m-t)^2$ has real
roots $t$.  Let $t_{{\rm min}}$ be the lesser root, and $t_{{\rm
max}}$ the greater root.  Then $S_\varepsilon \subset [0,t_{{\rm
max}}]$ if $t_{{\rm max}} \ge 0$, and $S_\varepsilon = \emptyset$ if
$t_{{\rm max}} < 0$.  By \lemref{bound}, we have $t_{{\rm max}}=0$ iff
$\ell^2+m^2=n^2/2$ and $m < n/2$, and $t_{{\rm max}} < 0$ iff
$\ell^2+m^2 > n^2/2$ and $m < n/2$.  In both of these cases, we have
${\mathcal P}(\ell/n,m/n)=0$ and $\Plt_\varepsilon = O(n^{-1})$.
Thus, we need only deal with the case $t_{{\rm max}} > 0$.

Suppose $t_{{\rm max}} > 0$ and $t_{{\rm min}} < 0$, i.e., $\ell^2+m^2
< n^2/2$.  It follows from \lemref{estim} that $k_{{\rm min}} =
O(n\varepsilon)$, and $k_{{\rm max}}=t_{{\rm max}} + O(n\varepsilon)$.
We will approximate the integral in (\ref{intap}) by
\begin{equation}
\label{int2}
\frac{1}{2}\int_{0}^{t_{{\rm max}}} {\frac{2}{\pi
\sqrt{(n-k)^2-2\ell^2-2(m-k)^2}}}\,dk.
\end{equation}
This approximation introduces further error.  To see how large the
error is, first rescale by a factor of $n$, so that the function under
the square root sign becomes $(1-\kappa)^2-2x^2-2(y-\kappa)^2$.
Around a root $r$, this function can be expanded as $-(\kappa-r)^2 \pm
2\sqrt{2}\sqrt{(y-1)^2-x^2}(\kappa-r)$ (with the sign depending on
which root $r$ is).  Because $|x|+|y|\le 1-\delta$, the coefficient of
$\kappa-r$ cannot become arbitrarily small.  Thus, for small
$\varepsilon$, the error introduced by the approximation is bounded by
a constant (depending on $\delta$) times
$$
\int_0^{\varepsilon} \frac{d\varepsilon'}{\sqrt{\varepsilon'}},
$$
and hence by $O(\varepsilon^{1/2}).$

One can evaluate the integral (\ref{int2}) explicitly, because
\begin{equation}
\label{inteval}
\int\!\frac{dk}{\sqrt{(n-k)^2-2\ell^2-2(m-k)^2}}
=
\tan^{-1}\!
\left(\frac{k+n-2m}{\sqrt{(n-k)^2-2\ell^2-2(m-k)^2}}\right)\!.
\end{equation}
As $k \rightarrow t_{{\rm max}}$, the right hand side of
(\ref{inteval}) approaches $\frac\pi 2$ (since the numerator of the
fraction is positive as its denominator approaches $0$).  We see that
(\ref{int2}) evaluates to
$$
\frac{1}{2}+
\frac{1}{\pi}\tan^{-1}
\left(\frac{2m-n}{\sqrt{n^2-2\ell^2-2m^2}}\right).
$$

The case with $t_{{\rm min}} \ge 0$ (i.e., $\ell^2+m^2\ge n^2/2$ and
$m > n/2$) is completely analogous, except the integral is over the
interval $[t_{{\rm min}},t_{{\rm max}}]$, rather than $[0,t_{{\rm
max}}]$.  This integral is $1$, so we get that $\Plt_\varepsilon = 1 +
O(n^{-1})+O(\varepsilon^{1/2})$, which agrees with ${\mathcal
P}(\ell/n,m/n)=1$.  This proves the desired result.
\end{proof}

We still need to prove that $\Plt_\varepsilon$ approximates
$\Pl(\ell,m;n+1)$.  We do that as follows:

\begin{thm}
\label{maininside}
Let $\delta > 0$ and $\varepsilon > 0$, such that $\varepsilon$ is
sufficiently small compared to $\delta$.  If $|\ell|+|m| \le
(1-\delta)n$ and $\ell+m \equiv n \pmod{2}$, then
$$
\Pl(\ell,m;n+1) = {\mathcal P}(\ell/n,m/n) 
+ O_\delta(\varepsilon^{1/2})
+ O_{\varepsilon,\delta}(n^{-1}).
$$
\end{thm}

\begin{proof}
We need to show that $\Plt_\varepsilon$ approximates
$\Pl(\ell,m;n+1)$.  Since \propref{gessel} implies that the creation
rates are all non-negative, and $\Plt_\varepsilon$ is the sum of a
subset of the creation rates appearing in the sum giving
$\Pl(\ell,m;n+1)$, the placement probability must be at least
$\Plt_\varepsilon$.

Also, given any point in the Aztec diamond, the north-going placement
probabilities at the four points obtained by rotating it by multiples
of $90^\circ$ about the origin sum to $1$.  This is true because by
rotational symmetry these placement probabilities are equal to the
placement probabilities in each of the four directions at the original
point, which must sum to 1.  This is the content of (\ref{rotation}),
except here it is expressed in terms of the placement probabilities,
rather than the asymptotic formula.

One can check by direct computation that ${\mathcal P}(x,y)+{\mathcal
P}(y,-x)+{\mathcal P}(-x,-y)+{\mathcal P}(-y,x)=1.$ If the difference
between $\Plt_\varepsilon$ and the placement probability were not $
O(\varepsilon^{1/2}) + O(n^{-1})$, then the four placement
probabilities would have to sum to more than $1$, which is impossible.
\end{proof}

The statement of \thmref{maininside} implies that away from the edges
of the diamond, the placement probabilities converge uniformly.
(Given any $\varepsilon > 0$, the theorem implies that if $n$ is
sufficiently large, then the placement probabilities are within a
constant multiple of $\varepsilon^{1/2}$ of their limiting values.  In
fact, the slightly awkward theorem statement is equivalent to
asserting uniform convergence; we state it that way because it seems
to be the form in which it is most naturally proved, given our setup.)
Thus, assuming \propref{expvanish}, we have very nearly proved
\thmref{main}.  In Section~\ref{sec-conclusion}, we will complete the
proof, using the following proposition:

\begin{prop}
\label{arcticbound}
For each $\varepsilon > 0$, there exists a positive constant $r < 1$
such that whenever $\ell^2 + m^2 > (1+\varepsilon)n^2/2$,
$$
\Pl(\ell,m;n+1) =
\begin{cases}
O(r^n)&\hbox{if $m < n/2$, and}\\
1+O(r^n)&\hbox{if $m > n/2$.}\\
\end{cases}
$$
\end{prop}

\begin{proof}
First suppose that $m < n/2$.  The desired result will follow from the
equation
\begin{equation}
\label{probsum}
\Pl(\ell,m;n+1) = \frac{1}{2}\sum_{k \ge 0}{\Cr(\ell,m-k;n+1-k)},
\end{equation}
together with the estimate given by \propref{expsmall}.  First, we
show that \propref{expsmall} applies to the creation rates appearing
in the sum.  Consider
\begin{equation}
\label{simplefn}
\frac{\ell^2+(m-k)^2}{(n-k)^2}
\end{equation}
as a function of $k$.  Its first derivative at $0$ is
$$
2\frac{\ell^2+m^2-mn}{n^3},
$$
which is greater than $0$ since $\ell^2+m^2 > n^2/2 > mn$.  The only
root of the derivative is
$$
\frac{\ell^2+m^2-mn}{m-n} < 0.
$$
Thus, the function (\ref{simplefn}) is increasing for $0 \le k < n$.
(Note that in (\ref{probsum}) we need only sum up to $k=(m+n)/2$,
since beyond that point $m-k<-(n-k)$ and hence
$\Cr(\ell,m-k;n+1-k)=0$.  Thus, $k$ never reaches the pole in
(\ref{simplefn}) at $n$.)  Therefore, $\ell^2+(m-k)^2 >
(1+\varepsilon)(n-k)^2/2,$ and \propref{expsmall} applies to bound the
creation rates in (\ref{probsum}).

Thus, for some constant $s$ between $0$ and $1$,
$$
\Pl(\ell,m;n+1) \le \sum_{k=0}^{(m+n)/2} O(s^{n-k}).
$$
This geometric series is bounded by $O(s^{(n-m)/2}) = O(s^{n/4}).$
This proves the desired bound, with $r=s^{1/4}$.

For $m > n/2$, we use the trick of summing the placement probabilities
at the four points obtained by rotating by multiples of $90^\circ$
about the origin.  As in the proof of \thmref{maininside}, the sum
must be $1$, and we know that three of the terms are $O(r^n)$.
Therefore, the fourth must be $1+O(r^n)$, as desired.
\end{proof}

\section{Exponential Sums}
\label{sec-exponential}

In the proof of \propref{firstapprox}, we needed to show that
$\Plt_\varepsilon$ is within $O(n^{-1})$ of the sum (\ref{firstsum});
to do so, we made use of an estimate whose proof was deferred
(\propref{expvanish}).  In this section, we will derive that estimate.
We begin with the following lemma.

\begin{lem}
\label{algroots}
Let $F(x_1,\dots,x_{n+1})$ be an algebraic function of $n+1$ variables
(defined on a subset of $\C^{n+1}$ to be specified shortly), and let
$S$ be a subset of $\C^n$.  Suppose that for each $(y_1,\dots,y_n) \in
S$, there exists an open set $U \subset \C$ such that as a function of
$x_{n+1}$, $F(y_1,\dots,y_n,x_{n+1})$ is (defined and) holomorphic on
$U$.  Then there is a constant $N$ such that for any $(y_1,\dots,y_n)
\in S$, if we regard $F(y_1,\dots,y_n,x_{n+1})$ as a function of
$x_{n+1}$ on the corresponding $U$, then it has at most $N$ roots in
$U$ if it is not identically zero.
\end{lem}

\begin{proof}
Since $F(x_1,\dots,x_{n+1})$ is algebraic, it satisfies an equation
\begin{equation}
\label{algeq}
\sum_{i=0}^{d}{p_i(x_1,\dots,x_{n+1})X^i} = 0,
\end{equation}
with $p_0,\dots,p_d$ polynomials (not all identically zero).  Let $N =
\max_i \deg p_i$.  We will show that $N$ has the desired property,
using induction on $n$.

We can choose the coefficients $p_i$ so that they have no
(non-constant) common factor.  Fix $y_1 \in \C$, and let $S' =
\{(y_2,\dots,y_n) \in \C^{n-1} : (y_1,\dots,y_n) \in S\}$.  Define
$G(x_2,\dots,x_{n+1}) = F(y_1,x_2,\dots,x_{n+1})$.  Since the
coefficients were taken to have no common factor, they do not all
vanish when we set $x_1=y_1$.  Their degrees do not increase when we
set $x_1=y_1$ (or when we remove common factors), so our lemma follows
by induction on $n$ (applied to $G$ and $S'$), assuming we can prove
it in the case $n=0$.

Suppose $n=0$.  Assuming $F$ is not identically zero, we can divide
(\ref{algeq}) by some power of $X$ to get an equation satisfied by $F$
with non-zero constant term, say $p_h(x_1)$.  (A priori, $F$ will
satisfy the new equation only where $F$ is non-zero.  However, since
$F$ is holomorphic on $U$, its zeros are isolated.  By continuity, it
satisfies the equation at its zeros as well as elsewhere.)  Then any
root of $F$ is a root of $p_h$, so $F$ has at most $\deg p_h$ roots,
and hence at most $N$ roots.
\end{proof}

\begin{lem}
\label{expbound}
The exponential sums
$$
\sum_{k \in I} \exp(2i\,\Phi(\ell,m-k;n-k))
$$
remain bounded (uniformly in $I$) as $n$ goes to infinity, where $I$
can be any subinterval of $S_\varepsilon$, as long as $|\ell|+|m| \le
(1-\delta)n$ for some fixed $\delta > 0$, and $\varepsilon$ is small
enough compared to $\delta$.
\end{lem}

\begin{proof}
To prove this, we will apply the Kusmin-Landau Theorem
(\thmref{kusmin-landau}).  \lemref{algebraic} says that
$d^2\Phi(\ell,m-k;n-k)/dk^2$ satisfies the conditions of
\lemref{algroots}, so there is an absolute upper bound for the number
of roots that it can have as a function of $k$ while $n$, $\ell$, and
$m$ are held fixed (unless it is identically zero for those values of
$n$, $\ell$, and $m$).  Before we apply the Kusmin-Landau Theorem, we
break $S_\varepsilon$ up into a bounded number of subintervals on
which $d\Phi(\ell,m-k;n-k)/dk$ is monotonic.

We have to look at the behavior of $d\Phi(\ell,m-k;n-k)/dk$.  As in
\lemref{phase}, set $k={\kappa}n$, $\ell=xn$, and $m=yn$.
\lemref{phase} says that as $n$ goes to infinity,
$d\Phi(\ell,m-k;n-k)/dk$ equals
\begin{equation}
\label{thetaform}
\theta\left({-x+y-{\kappa}},
{\sqrt {(1-{\kappa})^2-2\,(x^2+(y-{\kappa})^2)}
}\right)-\cr
\theta\left({1-{\kappa}-2\,x},
{\sqrt {(1-{\kappa})^2-2\,(x^2+(y-{\kappa})^2)}
}\right)
+O\left(\frac{1}{n}\right).
\end{equation}
We would like to show that when divided by $\pi$, (\ref{thetaform})
stays away from integers.

After (\ref{thetaform}) is divided by $\pi$, the only possible
integral values it can take on are $0$, $\pm 1$, and $\pm 2$ (assuming
$n$ is large enough).  If we ignore the $O(1/n)$ term, the rest of the
formula is the difference of the arguments of two points on the same
horizontal line (divided by $\pi$).  Thus, it cannot be $\pm 2$.  It
can be $0$ only if the points coincide or are on the horizontal axis.
It can be $\pm 1$ only if the points are on the horizontal axis.  The
points coincide iff $x+y=1$, which is impossible (since $|\ell|+|m|
\le (1-\delta)n$).  They are on the horizontal axis iff
\begin{equation}
\label{horaxis}
x^2+(y-{\kappa})^2=(1-{\kappa})^2/2.
\end{equation}
The definition (\ref{sdef}) of $S_\varepsilon$ implies that
$$
x^2+(y-{\kappa})^2 \le (1-\varepsilon)(1-{\kappa})^2/2,
$$
so no $k \in S_\varepsilon$ gives a ${\kappa}$ satisfying
(\ref{horaxis}).  (Note that $\kappa = 1$ is impossible since then
$|x|+|y| = |0|+|1| > 1-\delta$.)

In fact, the above argument, combined with continuity considerations,
shows that the two points cannot get arbitrarily close to each other
or the horizontal axis, and they clearly cannot get arbitrarily far
{}from the origin.  Thus, even taking into account the $O(1/n)$ term,
$(d\Phi(\ell,m-k;n-k)/dk)/\pi$ really does stay slightly away from
integers as $n \rightarrow \infty$.  Hence, the Kusmin-Landau Theorem
tells us that the exponential sums are bounded (uniformly in $I$).
\end{proof}

\begin{prop}
\label{expvanish}
The sum
$$
\sum_{k \in S_\varepsilon} {\frac{1}{\pi
\sqrt{(n-k)^2-2\ell^2-2(m-k)^2}}} \exp(2i\,\Phi(\ell,m-k;n-k)) 
$$
is $O(n^{-1})$ as $n$ goes to infinity, as long as $|\ell|+|m| \le
(1-\delta)n$ for some fixed $\delta > 0$, and $\varepsilon$ is small
enough compared to $\delta$.
\end{prop}

\begin{proof}
Let
$$
a(k) = {\frac{1}{\pi \sqrt{(n-k)^2-2\ell^2-2(m-k)^2}}}
$$
and
$$
b(k) = \sum_{k'=k_{{\rm min}}}^{k-1}\exp(2i\,\Phi(\ell,m-k';n-k')).
$$
For $k \in S_\varepsilon$, $a(k) = O(n^{-1})$ (by \lemref{maxbound})
and $b(k)$ is bounded (by \lemref{expbound}).  Suppose $|b(k)| \le B$
for all $k \in S_\varepsilon$.

To bound the sum in the statement of the proposition, we will apply
summation by parts.  We have
\begin{eqnarray*}
\sum_{k \in S_\varepsilon} {\frac{\exp(2i\,\Phi(\ell,m-k;n-k))}{\pi
\sqrt{(n-k)^2-2\ell^2-2(m-k)^2}}} 
& = & \sum_{k \in S_\varepsilon}{a(k)(b(k+1)-b(k))} \\
& = & \sum_{k \in S_\varepsilon}a(k)b(k+1) - \sum_{k \in
S_\varepsilon}a(k)b(k) \\ 
& = & \sum_{k=k_{{\rm min}}+1}^{k_{{\rm max}}+1}\!\!a(k-1)b(k) -
\sum_{k=k_{{\rm min}}}^{k_{{\rm max}}}\!a(k)b(k) \\
& = & \sum_{k=k_{{\rm min}}+1}^{k_{{\rm max}}}\!\!b(k)(a(k-1)-a(k)) +
O(n^{-1}).
\end{eqnarray*}

This sum is bounded in absolute value by $B \sum_{k}{|a(k-1)-a(k)|} +
O(n^{-1})$.  The function $a(k)$ is monotonic on $(-\infty,2m-n)$ and
$(2m-n,\infty)$ (on the subintervals where it is real, of course), so
within each of these intervals, the sum $\sum_{k}{|a(k-1)-a(k)|}$
telescopes.  The boundary terms are $O(n^{-1})$, and hence the entire
sum is $O(n^{-1})$.
\end{proof}

\section{Conclusion of the Proof}
\label{sec-conclusion}

The results proved in the preceding three sections give us
\thmref{maininside}, a weakened version of Theorem~\ref{main}, in
which we are restricted to estimating the placement probabilities at
normalized locations $(x,y)$ with $|x|+|y|\le 1-\delta$ for some fixed
$\delta>0$.  That is, we are required to keep $(x,y)$ from getting too
close to the boundary of the diamond.  Here we will show how the
restriction on $(x,y)$ can be relaxed, provided that we are careful to
stay away from the points $(\pm\frac12, \frac12)$.

Fix $\delta>0$, and consider the region in the Aztec diamond of order
$n$ defined (relative to normalized coordinates) by the constraint
$x^2+y^2>\frac12+\delta$.  This region splits up into four pieces.
\propref{arcticbound} tells us that the north-going placement
probabilities tend uniformly to 1 in the northern piece and to 0 in
the other three pieces.  The only regions that are not covered by this
method are four small curvilinear trapezoids near the points $(\pm
\frac12, \pm \frac12)$, defined by the inequalities
$x^2+y^2\le\frac12+\delta$ and $1-\delta < |x|+|y| \le 1$.  If $(x,y)$
stays away from these four points as $n$ goes to infinity, then we can
indeed conclude that the placement probabilities for north-going
dominos at location $(x,y)$ are as claimed in Theorem~\ref{main}.
This completes the proof of the main theorem, except near the points
$(\pm\frac12,-\frac12)$, which we will now deal with.

Let $R$ be the subregion of the Aztec diamond of order $n$ consisting
of the two lower of the four curvilinear trapezoids defined by
$x^2+y^2 \le \frac12 + \delta$ and $1-\delta < |x|+|y| \le 1$.  In
$R$, we use the inequality (\ref{incbound}).  It says that for $h \ge
0$,
$$
\Pl(\ell,m;n) \le \Pl(\ell,m+h;n+h).
$$
If $n$ is sufficiently large, then for each point $(\ell,m)$ in $R$,
there exists an $h$ such that the point $(\ell,m+h)$ of the diamond of
order $n+h$ has normalized coordinates satisfying $x^2+y^2 \le \frac12
+ \delta$, $1-2\delta < |x|+|y| \le 1-\delta$, and $y < 0$.  Let $S$
be the set of all $(x,y)$ satisfying these three constraints.
Inequality (\ref{incbound}) tells us that the placement probabilities
within $R$ are at most as large as those within $S$.  However, the
part of \thmref{main} that we have already proved gives estimates for
the placement probabilities in $S$, and shows that they tend uniformly
to $0$ as $\delta \rightarrow 0$.  (To see the convergence to $0$ most
easily, look at the level curves of the placement probabilities.)  We
thus conclude that as $\delta \rightarrow 0$, the placement
probabilities in $R$ tend uniformly to $0$.  This completes the proof
of \thmref{main}.

Unfortunately, our techniques do not give us an explicit bound for the
difference between the placement probabilities and the arctangent
formula in an Aztec diamond of a given order.  This is not because the
methods are inherently ineffective; rather, it is because we have not
determined the dependence on $\varepsilon$ in the $O(n^{-1})$ term of
the error bound in \thmref{maininside}.  To determine it, we would
have to do so for the error term in \propref{crest}, which seems more
trouble than it would be worth (but could perhaps be done).

Using these techniques, we can also prove the arctic circle theorem of
\cite{circle}.  One direction, that the regions outside the inscribed
circle are indeed frozen, follows from \propref{arcticbound}.  To see
this, consider the region $R$ defined (relative to normalized
coordinates) by $x^2+y^2>\frac{1}{2}+\varepsilon$, with $\varepsilon >
0$.  The number of domino spaces in this region of an Aztec diamond of
order $n$ is less than $n^2$, so the probability that a
non-north-going domino will appear in the subregion with $y >
\frac{1}{2}$, or that a north-going domino will appear in the rest of
$R$, is exponentially small, by \propref{arcticbound}.  From this,
we see that with probability approaching $1$ (as $n$ goes to
infinity), all the dominos in $R$ will be aligned in brickwork
patterns, and thus contained in the polar regions.  This is half of
the arctic circle theorem.

The other direction, that the polar regions almost never extend
substantially into the interior of the inscribed circle, requires an
additional result for its proof.  Intuitively, it follows from our
main theorem, which tells us that inside the inscribed circle all four
types of placement probabilities are positive.  This trivially implies
that the polar regions cannot almost always cover a given part of the
interior of the circle, but showing that they almost never do is
harder.  We will prove it in subsection~\ref{ssec-ACT}.

\section{Consequences of the Theorem}
\label{sec-consequences}

\subsection{Height functions}
\label{ssec-height}

Height functions for domino tilings, which were introduced in the
mathematics literature in \cite{conway} (and independently in a
slightly different form in the physics literature in \cite{levitov}),
are a very useful device in the study of tilings of simply-connected
subsets of the plane.  (A more general approach to height functions
can be found in \cite{spaces}.)  In any such region that can be tiled
by dominos, the number of enclosed white squares and the number of
enclosed black squares under an alternating coloring of the squares
must clearly be equal.  It follows that if one travels around the
boundary of the region counterclockwise, then one will see a black
square on one's left half the time and a white square on one's left
half the time; to see why, notice that the edges of the square grid
that lie within the region pair sides of black squares with sides of
white squares, so the excess of unpaired (i.e., boundary) sides of
black squares over unpaired sides of white squares is four times the
excess of black squares over white squares.  As one travels around the
boundary, any temporary excess of one kind of square over the other
that is observed along the way represents a ``debt'' that will
eventually have to be paid.  Moreover, the same is true simultaneously
for all the boundaries of all the simply-connected regions that are
formed by suitable subsets of the tiles in question.  Height functions
provide a uniform framework for keeping track of all these debts
simultaneously.

If $R$ denotes a finite, simply-connected region composed of lattice
squares that have been alternately colored black and white, a {\sl
height function} on $R$ is an integer-valued function $h$ on the
vertices of the lattice squares which satisfies the following two
properties for adjacent vertices $u$ and $v$: first, if the edge from
$u$ to $v$ is part of the boundary of $R$, then $|h(u)-h(v)|=1$, and
second, if the edge from $u$ to $v$ has a black square on its left,
then $h(v)$ is either $h(u)+1$ or $h(u)-3$.  It is not hard to show
that such a function necessarily satisfies a discrete Lipschitz
condition: if vertices $u$ and $v$ are at distance $d$ from each other
in the sup-norm, then $|h(u)-h(v)| \leq 2d+1$.  Note also that if
$h(\cdot)$ is a height function, then so is $h(\cdot)+C$ for any
integer $C$.

\begin{figure}
\begin{center}
\PSbox{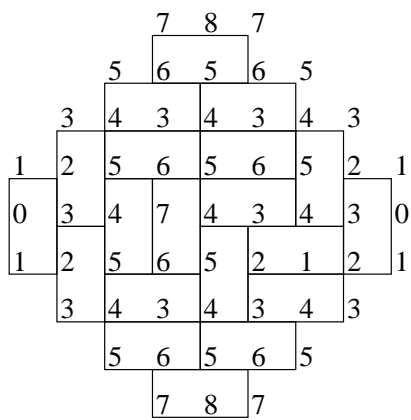}{2.0in}{2.25in}
\end{center}
\caption{A height function for a domino tiling.}
\label{fig-typical}
\end{figure}

Every height function on $R$ determines a domino tiling of $R$,
consisting of dominos that occupy all the domino spaces that are
bisected by edges $uv$ with the property that $|h(u)-h(v)|=3$.
Conversely, every domino tiling of $R$ arises in this way from a
height function on $R$ that is unique modulo addition of a global
constant.  We can remove this ambiguity by constraining a particular
vertex on the boundary of $R$ to have some particular integer value as
its height; then every domino tiling of $R$ has a unique height
function subject to this constraint, and what is more, all these
height functions agree with one another on the boundary of $R$.  For
instance, in the case of the Aztec diamond of order $n$, we set things
up so that the middle vertex on the west edge of the diamond has
height 0 and the middle vertex on the northern edge has height $2n$.
(Note that this differs by 1 from the height function convention used
in \cite{alternating}.)  Then the heights for a typical domino tiling
of an Aztec diamond are as shown in Figure~\ref{fig-typical}.  (The
shading convention for the lattice squares is the same as that in
Section~\ref{sec-introduction}, i.e., so that the leftmost square of
each row in the top half of the diamond is white.)

One can also develop an analogous theory of height functions for other
sorts of tilings, for example, tilings of regions in the triangular
lattice by lozenges (two unit equilateral triangles joined along an
edge).  This theory is simpler geometrically than that for domino
tilings; for the details, see \cite{conway}.  (See also \cite{bh} for
an independent, earlier development of height functions for this
lattice in physics.)  Height functions can furthermore be applied to
the square ice model studied by Lieb, as is shown in \cite{bcsos}.
The results of subsections~\ref{ssec-robustness} and
\ref{ssec-variance} generalize straightforwardly to other sorts of
height functions.  However, because the focus of this article is on
domino tilings, we will not go into the details of the generalization.

Suppose $u$, $v$, and $w$ are three consecutive vertices along a path
in a simply-connected region $R$ that is tiled by dominos, such that
neither the edge $uv$ nor the edge $vw$ bisects a domino.  Then
$|h(u)-h(v)|=|h(v)-h(w)|=1$, with $h(w)=h(u)$ if the three points are
collinear and $h(w)=h(u) \pm 2$ if they are not.  This principle makes
it fairly easy to go through the tiling, assigning heights to the
vertices.  Alternatively, one can use this method just to find the
heights along the boundary, and then find the heights in the interior
by the following procedure.  To determine the height of a particular
vertex in the interior of a tiled region, start at the point on the
boundary of the region due north of the vertex (whose height is
independent of the tiling) and proceed downward, subject to the
following rule: when one travels southward along an edge that bisects
a north-going (resp.\ south-going) domino in the tiling, the height
decreases (resp.\ increases) by 3, whereas, when one travels southward
along an edge that bisects a north-going (resp.\ south-going) domino
space that is not occupied by a domino in the tiling, the tiling, the
height increases (resp.\ decreases) by 1.  A similar rule can be
formulated for describing how the height changes as one travels
horizontally through the interior of the diamond.  The fact that these
rules are consistent with each other is a consequence of the fact that
any region that can be tiled by dominos must contain exactly equal
numbers of black and white squares.

If we take the average of all the (finitely many) height functions
associated with the different tilings of a region, we get a
real-valued function on the vertices called the {\sl average height
function}.  As a consequence of the rule described in the preceding
paragraph, one can give a simple description of how the average height
changes from vertex to vertex, in terms of the placement probability
$p$ associated with the domino space that is bisected by the edge that
connects the two vertices.  For instance, if $u$ and $v$ are
neighbors, with $u$ to the north or west of $v$, then the average
height at $v$ is equal to the average height at $u$ plus $4p-1$ if
edge $uv$ bisects a south-going or west-going domino space, while the
average height at $v$ is equal to the average height at $u$ minus
$4p-1$ if edge $uv$ bisects a north-going or east-going domino space.

Here we are interested in the asymptotic behavior of the average
height function for domino tilings of the Aztec diamond.  The
arctangent formula tells us that these probabilities $p$ are slowly
varying, so the average height function is locally approximated by
functions of the form $ax+by+c$ (with $a$, $b$, and $c$ slowly
varying).  We call the pair $(a,b)$ the {\sl tilt} of the plane
$z=ax+by+c$.  Let us normalize our height functions by dividing
through by $n$, both in the domain and in the range.  Thus, in the
limit we expect to see some sort of function ${\mathcal
H}(\cdot,\cdot)$ on $\{(x,y): |x|+|y| \leq 1\}$ satisfying the
piecewise-linear boundary condition ${\mathcal H}(x,y) = 1-x^2+y^2$
for $|x|+|y|=1$, as well as a Lipschitz condition with constant 2
relative to the sup-norm distance.  In addition, the formulation in
the previous paragraph of how the average height changes when moving
between vertices tells us that we should have $\frac{\partial{\mathcal
H}}{\partial y}=2(p_n-p_s)$ and $\frac{\partial{\mathcal H}}{\partial
x}=2(p_w-p_e)$, where $p_n$, $p_s$, $p_e$, and $p_w$ are the
north-going, south-going, east-going, and west-going placement
probabilities at $(x,y)$, respectively.

It can be shown (although we do not prove this here) that the domino
shuffling algorithm of \cite{alternating} can be interpreted directly
in terms of height functions, and that half of the values of the
average height function for the diamond of order $n+1$ are equal to
certain corresponding values of the average height function for the
diamond of order $n$.  Hence the average height functions for the
Aztec diamonds of orders $n$ and $n+1$ cannot be too far apart.
However, such considerations are not sufficient to yield a proof that
the normalized average height functions converge to a continuum limit.

The arctangent formula gives us the strength we need in order to
conclude that a limit exists.  Recall that the average height function
can be derived by taking cumulative sums and differences of local
placement probabilities, with various coefficients.  Taking this
assertion to the limit as $n \rightarrow \infty$, and using the fact
that the placement probabilities approach a continuum limit, we see
that the normalized average height function also approaches a limit.
(It is true that the errors in the placement probabilities are going
to add, and that there are more and more of them as $n$ gets large,
but each individual error is going to be small, so that when we divide
by $n$ the normalized error is small as well.)  The limit must be some
function ${\mathcal H}(x,y)$ (defined for $|x|+|y| \leq 1$) with the
property that $\frac{\partial {\mathcal H}}{\partial x}=2{\mathcal
P}(y,-x)-2{\mathcal P}(-y,x)$ and $\frac{\partial {\mathcal
H}}{\partial y}=2{\mathcal P}(x,y)-2{\mathcal P}(-x,-y)$ for all
$(x,y)$ in the interior of its domain.  (Here, we have expressed the
placement probabilities near $(x,y)$ for all four domino orientations
in terms of ${\mathcal P}$ via rotational symmetry.)  That is, the
tilt of the tangent plane at a point (associated with the average
height function) can be expressed in terms of the local placement
probabilities for random domino tilings.  This means that we ought to
be able to reconstruct the function ${\mathcal H}(\cdot,\cdot)$ from
Theorem 1 via integration, making use of the known boundary conditions
satisfied by ${\mathcal H}$.  If we do this, it turns out that
${\mathcal H}$ can be written in closed form, and indeed, a formula
for ${\mathcal H}$ can be written especially compactly if one makes
use of the formula for ${\mathcal P}(\cdot,\cdot)$ itself.
Specifically, one can verify that the following formula for ${\mathcal
H}(\cdot,\cdot)$ holds:

\begin{prop}
\label{heightformula}
The normalized average height functions for large Aztec diamonds
converge uniformly to
$$
{\mathcal H}(x,y) = 2\left(y{\mathcal P}(x,y)-y{\mathcal
P}(-x,-y)+(1-x){\mathcal P}(-y,x)+(1+x){\mathcal P}(y,-x)\right).
$$
\end{prop}

\begin{proof}
We simply check that this formula satisfies the differential equations
and boundary conditions.
\end{proof}

Within the temperate zone, the average height function is real
analytic; in each of the polar regions, it is an affine function of
$x$ and $y$.  On the arctic circle itself, away from the points $(\pm
\frac12, \pm \frac12)$, the function ${\mathcal H}(x,y)$ is
differentiable but not twice-differentiable.  It takes the value 2 at
the points $(0, \pm 1)$ and the value 0 at the points $(\pm 1,0)$,
with piecewise-linear behavior on the boundary of the normalized
diamond.  The level set ${\mathcal H}=1$ consists of the two line
segments joining midpoints of opposite sides of the normalized
diamond.

One may ask, for diamonds of finite order $n$, how closely the
distribution on height functions is clustered around its mean value.
We will see in the next subsection that the standard deviation of the
height at any fixed location in the Aztec diamond of order $n$ is at
most $8\sqrt{n}$.  However, numerical evidence suggests that, at the
center of the diamond, the standard deviation of the height is much
smaller---more like $\log{n}$, or perhaps even less than that.

Our formulas for $\frac{\partial {\mathcal H}}{\partial x}$ and
$\frac{\partial {\mathcal H}}{\partial y}$, in combination with the
arctangent formula, yield (within the temperate zone) the equation
\[
\frac{\partial^2 {\mathcal H}}{\partial y^2}
-\frac{\partial^2 {\mathcal H}}{\partial x^2}
= \frac{8}{\pi \sqrt{1-2x^2-2y^2}}.
\]
We can rewrite this equation in a slightly more illuminating way.  For
$t>0$ and $|x|+|y|\le t$, define
$${\overline{\mathcal H}}(x,y,t)=t{\mathcal H}(x/t,y/t).$$
That is, we undo the scaling introduced in Section~\ref{sec-introduction}.
Then off the arctic circle we have
$$
\frac{\partial^2 {\overline{\mathcal H}}}{\partial y^2}-
\frac{\partial^2 {\overline{\mathcal H}}}{\partial x^2}
= 8u(x,y,t),
$$
where
$$
u(x,y,t) = 
\begin{cases}
\frac{1}{\pi\sqrt{t^2-2x^2-2y^2}}&\hbox{if $x^2+y^2<t^2/2$, and}\\
0&\hbox{if $x^2+y^2\ge t^2/2$.}\\
\end{cases}
$$
This function is a fundamental solution to the wave equation in two
dimensions, with speed of propagation $1/\sqrt{2}$.  That is, $u$ is a
distribution satisfying
$$
\frac{\partial ^2 u}{\partial t^2} = \frac{1}{2}\left(
\frac{\partial^2 u}{\partial x^2} + \frac{\partial^2 u}{\partial y^2}
\right),
$$
$$
u(x,y,0) = 0,
$$
and
$$
\frac{\partial u}{\partial t}(x,y,0) = \delta(x,y),
$$
where $\delta$ is the (two-dimensional) Dirac delta function.  (For
more details on fundamental solutions to the wave equation, see
\cite[p.~164]{rauch}.)

Note that \propref{crest} shows that, except for an oscillating
factor, the creation rates also behave like $2u$.  William Jockusch
has shown in personal communication how to use a generating function
developed in \cite{gip} to explain this behavior, by viewing the
creation rates as numerical approximations to a solution of the wave
equation.

He has also pointed out that from his methods, one ought to be able to
deduce a weak version of Theorem~\ref{main}.  More specifically, one
should be able to show that in any macroscopic subregion of a randomly
tiled Aztec diamond of order $n$ (i.e., any subregion of size
comparable to that of the diamond), the expected density of
north-going dominos is within $o(1)$ of that predicted by integrating
the arctangent formula; in particular, this would suffice to prove
\propref{heightformula}.  Unfortunately, his methods would not rule
out the possibility of small-scale fluctuations in the placement
probabilities, such as one gets if one looks at placement
probabilities for all horizontal domino spaces rather than just the
north-going ones.

\subsection{Robustness}
\label{ssec-robustness}

\begin{figure}
\begin{center}
\PSbox{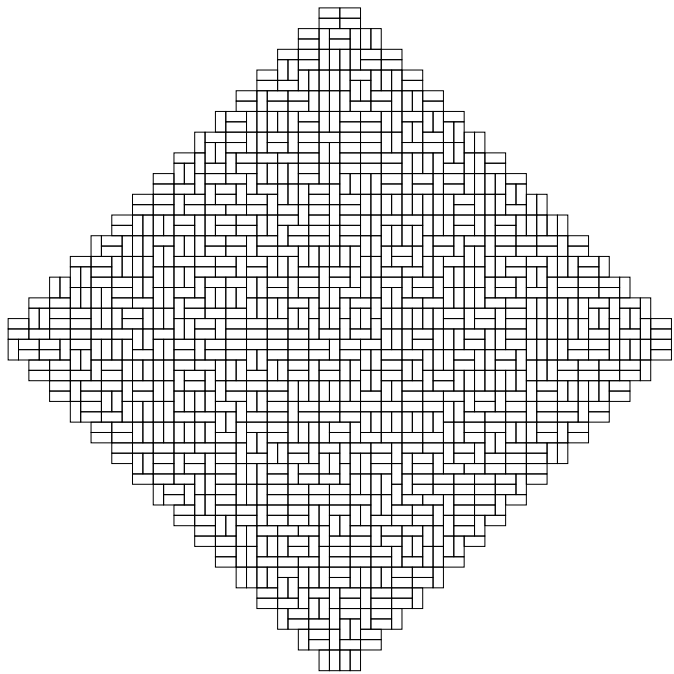}{250pt}{250pt}
\end{center}
\caption{A modified Aztec diamond with far more tilings.}
\label{fig-fake1}
\end{figure}

\begin{figure}
\begin{center}
\PSbox{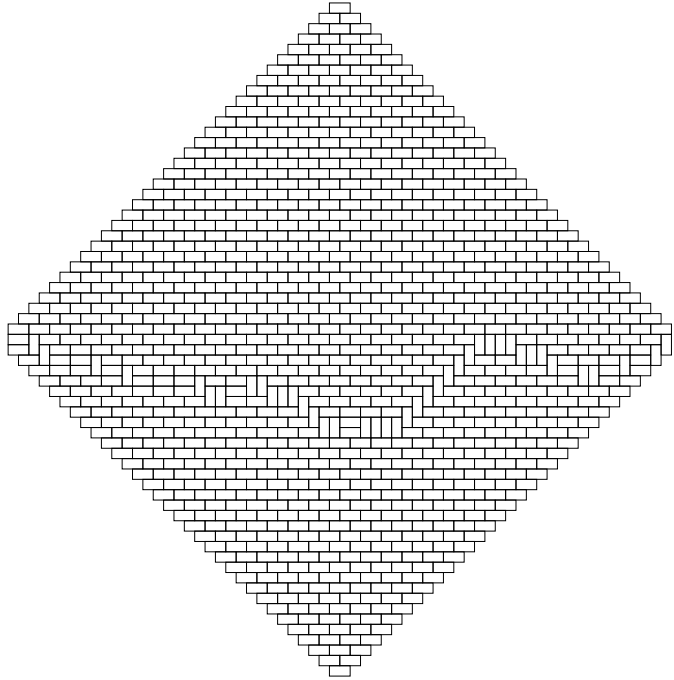}{250pt}{260pt}
\end{center}
\caption{A modified Aztec diamond with far fewer tilings.}
\label{fig-fake2}
\end{figure}

\begin{figure}
\begin{center}
\PSbox{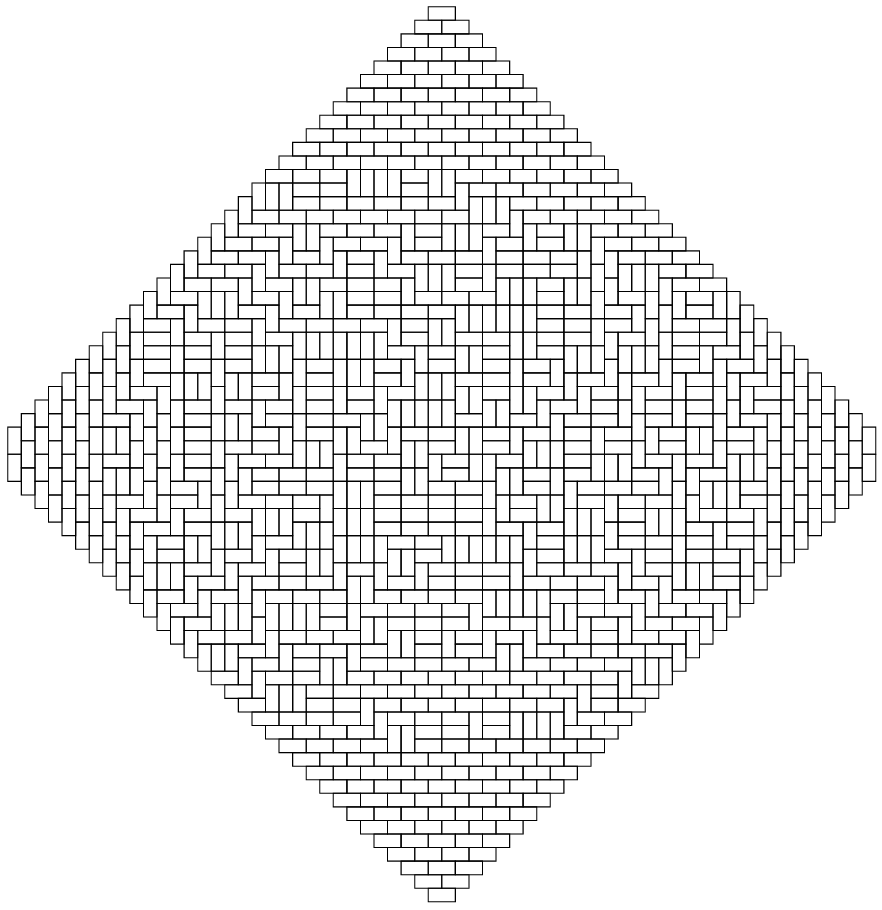}{250pt}{260pt}
\end{center}
\caption{A modified Aztec diamond with Aztec-like height function.}
\label{fig-real1}
\end{figure}

The formula for the average height function that was derived in the
preceding section from the arctangent formula applies not only to
Aztec diamonds, but also to all regions that approximate them in a
suitable sense.  (Here, as hereafter, the term ``region,'' without
qualifiers, should be understood to refer to finite regions in the
plane that are unions of lattice squares and can be tiled with
dominos.)  It is not enough that the region being tiled should have a
boundary that is roughly ``Aztec'' in shape.  For instance,
Figure~\ref{fig-fake1} shows a random domino tiling of a region
obtained from the Aztec diamond of order 32 by adding some squares
along its boundary, while the region that is shown in
Figure~\ref{fig-fake2}, also studied in \cite{remark}, was obtained by
adding an extra row of length 64 in the middle of the diamond.  (These
random tilings were obtained via the method described in \cite{exact}
and are indeed truly random, to the extent that pseudo-random number
generators can be trusted.)  In neither case do we get behavior
consistent with the arctangent formula.  On the other hand,
Figure~\ref{fig-real1} shows a random tiling of an Aztec diamond to
which {\sl two} rows of length 64 have been added, and the resemblance
to Figure~\ref{fig-aztec} is evident.

The sense of mystery dissolves if one considers the behavior of the
height function along the boundary in each of the three cases.  In the
first case, the height function is nearly constant along the boundary;
in the second, the direction of change of the height function is the
same along the southwest and northwest edges (and also the same along
the southeast and northeast edges); and in the third, the direction of
change of the height function switches as one rounds any of the four
corners of the region.  Since it is the third situation that resembles
the boundary behavior of height functions of Aztec diamonds, it is not
surprising that the third situation should also give behavior in the
interior that is similar to what one sees for Aztec diamonds.  (In
fact, if we view the third region as an Aztec diamond of order 33 with
two vertical dominos removed, then since almost all tilings of the
Aztec diamond contain those two dominos, the local statistics in the
third region differ little from those in the Aztec diamond of order
$33$.)

Note, incidentally, that if instead of adding a row of length 64, as
we did in Figure~\ref{fig-fake2}, we removed a row of length 64, then
the resulting region is easily seen to have a unique tiling,
consisting entirely of horizontal brickwork.  Although this situation
may seem trivial, it can shed some light on what is happening in
Figure~\ref{fig-fake2}.  The horizontal brickwork pattern seen almost
everywhere in Figure~\ref{fig-fake2} is the unique arrangement of
dominos such that the height increases (or decreases, depending on
whether the dominos in the pattern are north-going or south-going) as
quickly as possible as one moves vertically.  In an Aztec diamond with
a row of length 64 removed, the heights on the boundary are such that
the only way to span the gap between the heights on the lower half of
the boundary and those on the upper half is to change at this rate.
Thus, the only way to tile the region is with a brickwork pattern.  In
the case of Figure~\ref{fig-fake2}, the occurrence of an extra row of
length 64 gives the height function a tiny bit of slack, and we can
see where this slack gets used by following the fault-line that runs
{}from left to right.

Of course, we could have predicted ahead of time that small
modifications of the shape of the boundary can have a drastic impact
on the tiling situation, since for instance adding a single square to
a region (or removing a single square) can create a region with odd
area, which cannot be tiled at all.  Hence, we will want to assume
that all the regions we discuss actually admit tilings by dominos, as
stipulated in the first paragraph of this subsection.

We will show in this subsection that regions similar to Aztec
diamonds, such as Figure~\ref{fig-real1}, have approximately the same
average height functions as the Aztec diamonds they resemble.  This
will follow as a consequence of a more general result, asserting that
the value of the average height function depends in a continuous
manner on the values of the fixed heights along the boundary.  That
is, if one modifies the shape of the boundary in such a way that the
height function along the new boundary, when plotted in three
dimensions (the two original dimensions plus a third dimension for
height), is close to the graph of the old, the average heights of
vertices in the interior will not change very much.

Before we can do this, we first prove a general monotonicity result
about height functions.  The idea for this approach was suggested by
Robin Pemantle in personal communication.  Let $R$ denote a
simply-connected region in the plane with some fixed checkerboard
coloring, and let $V$ be the set of vertices in $R$.  Let $V'$ be a
subset of $V$ that contains all the vertices on the boundary of $R$;
we will assume that $V'$ is connected, in the sense that the subgraph
of the square grid induced by the vertex set $V'$ is connected.  A
{\sl partial height function} is a function $f : V' \rightarrow \Z$
subject to the local constraint that if $u$ and $v$ are adjacent
vertices such that the directed edge from $u$ to $v$ has a black
square on its left, then $f(v)-f(u)$ is either $1$ or $-3$.  It is
called {\sl complete height function} if it is defined on all of $V$;
a complete height function $\hat{f}$ {\sl extends} a partial height
function $f$ if it agrees with $f$ where $f$ is defined.

Throughout this subsection (and the next), $H$ will denote a complete
height function chosen at random (according to some distribution);
thus, for any vertex $v$, $\Exp[H(v)]$ (the expected value of $H(v)$)
is the value at $v$ of the average height function.

Given a connected subset $V'$ of $V$ that contains all the boundary
vertices, and given a partial height function $f$ on $V'$, we let
$\mu_f$ denote the uniform distribution on the set of complete height
functions that extend $f$ to $V$.

\begin{lem}
\label{monotone}
If $f$ and $g$ are two partial height functions defined on $V'$ and
agreeing modulo $4$, with $f \leq g$, then $\mu_f$ is stochastically
dominated by $\mu_g$.  That is, there exists a probability measure
$\pi$ on the set of pairs $(\hat{f},\hat{g})$ of complete height
functions extending $f$ and $g$ respectively, such that
$$
\sum_{\hat{g}} \pi(\hat{f},\hat{g})=\mu_f(\hat{f}),
$$
$$
\sum_{\hat{f}} \pi(\hat{f},\hat{g})=\mu_g(\hat{g}),
$$
and
$$
\pi(\{(\hat{f},\hat{g}): \hat{f} \leq \hat{g}\}) = 1.
$$
\end{lem}

\begin{proof}
We use induction on the size of $V \setminus V'$ (holding $V$ fixed
and varying $V'$).  The case where this set is empty is trivial.
Assume that the lemma is true whenever $|V \setminus V'| = k-1$, and
suppose we have a situation in which $|V \setminus V'| = k$.  It
clearly suffices to consider the case in which $f(v) < g(v)$ for some
vertex $v$ in $V'$ that is adjacent to at least one vertex in $V
\setminus V'$.  Let $w$ be a vertex in $W = V \setminus V'$ adjacent
to $v$.

Given that $f(v)$ has some specific value, any extension of $f$ to
$V'' = V' \cup \{w\}$ would have to give $w$ height $h$ or $h-4$ (for
some particular $h$ whose value we don't care about---it's $f(v)$ plus
or minus 1 or 3), while any extension of $g$ to $V''$ would have to
give $w$ height $h'$ or $h'-4$ (with $h'$ determined from $g(v)$ the
same way $h$ is determined from $f(v)$).  Because $f$ and $g$ agree
modulo $4$ on $V'$ and $h' > h$, we have $h'-4 \ge h$.

Let $f'_1$ and $f'_2$ be the two extensions of $f$ to $V''$ that
assign $w$ height $h$ and $h-4$, respectively, and let $g'_1$ and
$g'_2$ be the two extensions of $g$ to $V''$ that assign $w$ height
$h'$ and $h'-4$, respectively.  (If such extensions do not exist, it
is not a problem, as we will see below.)  The distribution $\mu_f$ is
a weighted superposition of $\mu_{f'_1}$ and $\mu_{f'_2}$, where the
$i$th term ($i=1$ or $2$) is given weight proportional to the number
of extensions of $f'_i$ to $V$ (which should be taken to be zero in
the case where the extension to $V''$ does not exist).  Similarly,
$\mu_g$ is a superposition of $\mu_{g'_1}$ and $\mu_{g'_2}$.  Since
$f'_i \leq g'_j$ for all $i,j$ in $\{1,2\}$, and $h \equiv h'
\pmod{4}$, we can use our induction hypothesis to conclude that
$\mu_{f'_i}$ is stochastically dominated by $\mu_{g'_j}$ for all
$i,j$, which implies that $\mu_f$ is stochastically dominated by
$\mu_g$, as was to be shown.
\end{proof}

\begin{cor}
\label{concor}
If $f$ and $g$ are two partial height functions on $R$ defined on $V'$
and agreeing modulo $4$, with $f \leq g+4$, then for all $v$,
$\Exp[H(v)]$ under the measure $\mu_f$ is at most 4 more than
$\Exp[H(v)]$ under the measure $\mu_g$.
\end{cor}

\begin{proof}
Apply \lemref{monotone} to the partial height functions $f$ and $g+4$.
\end{proof}

For applications of \lemref{monotone} and \corref{concor}, it is
important to note that the values of height functions on connected
regions are determined modulo $4$, given the value at any one point,
because the defining conditions for a height function imply that if
two height functions agree modulo $4$ at any point, then they do so at
each neighboring point.  Also, given two partial height functions
defined on different sets, we say that they agree modulo $4$ if all
height functions extending them agree modulo $4$.

\begin{prop}
\label{robustness}
Suppose that $R_1$, $R_2$ are two simply-connected regions in the
plane, with mandated partial height functions $f_1$, $f_2$ along their
boundaries that agree modulo $4$, such that every vertex $v_1$ on the
boundary of $R_1$ is within sup-norm distance $\Delta_1$ of some
vertex $v_2$ on the boundary of $R_2$, and vice versa, and such that
whenever vertices $v_1$ and $v_2$ on the respective boundaries are
within sup-norm distance $\Delta_1$ of each other, the heights
$f_1(v_1)$ and $f_2(v_2)$ are within $\Delta_2$ of each other.  Then,
for any $v$ in $R_1 \cap R_2$, the expected value of $H(v)$ under
$\mu_{f_1}$ and the expected value of $H(v)$ under $\mu_{f_2}$ differ
by at most $2\Delta_1 + \Delta_2 + 1$.
\end{prop}

\begin{proof}
Let $f_{1,{{\rm max}}}$ be the highest extension of $f_1$ to $R_1$,
let $f_{1,{{\rm min}}}$ be the lowest extension of $f_1$ to $R_1$, and
define $f_{2,{{\rm max}}}$ and $f_{2,{{\rm min}}}$ similarly.  (It is
not hard to show that the complete height functions extending a given
partial height function form a lattice under the usual partial
ordering, so it makes sense to talk about the highest and lowest
extensions.)  Let $v$ be on the boundary of $R_1 \cap R_2$ (and hence
on the boundary of $R_1$ or $R_2$).  If $v$ is on the boundary of
$R_1$, then we can find a nearby $w$ on the boundary of $R_2$ so that
\begin{eqnarray*}
f_{1,{{\rm max}}}(v) & = & f_1 (v) \\
& \leq & f_2 (w) + \Delta_2 \\
& = & f_{2,{{\rm min}}}(w) + \Delta_2 \\
& \leq & f_{2,{{\rm min}}}(v) + 2 \Delta_1 + 1 + \Delta_2,
\end{eqnarray*}
while if $v$ is on the boundary of $R_2$, then we can find a nearby
$w$ on the boundary of $R_1$ so that
\begin{eqnarray*}
f_{1,{{\rm max}}}(v) & \leq & f_{1,{{\rm max}}}(w) + 2 \Delta_1 + 1 \\
& = & f_1(w) + 2 \Delta_1 + 1 \\
& \leq & f_2(v) + \Delta_2 + 2 \Delta_1 + 1 \\
& = & f_{2,{{\rm min}}}(v) + \Delta_2 + 2 \Delta_1 + 1.
\end{eqnarray*}
Since the two height functions agree modulo $4$ at $v$, $f_{1,{{\rm
max}}} (v) \leq f_{2,{{\rm min}}} (v) + 4K$, where $4K$ is the
greatest multiple of 4 that is less than or equal to $2\Delta_1 +
\Delta_2 + 1$.  It follows from this (and the corresponding inequality
in the other direction) that if $f_1'$ is any extension of $f_1$ to
$R_1$ and $f_2'$ any extension of $f_2$ to $R_2$, then for each $v$ on
the boundary of $R_1 \cap R_2$, $f_1'(v)$ differs from $f_2'(v)$ by at
most $4K$.

Now let $v$ be any vertex in $R_1 \cap R_2$.  If we compute
$\Exp[H(v)]$ by conditioning on the heights on the boundary of $R_1
\cap R_2$, then it follows from \corref{concor} that the expected
value of $H(v)$ under $\mu_{f_1}$ differs by at most $4K$ (and hence
at most $2\Delta_1 + \Delta_2 + 1$) from its expected value under
$\mu_{f_2}$.
\end{proof}

As an application of this result, we may consider a modification of
the Aztec diamond of order $n$, whose symmetric difference with the
true Aztec diamond of order $n$ is a narrow fringe around the border
of the true diamond, of width $o(n)$.  Suppose that the black and
white squares of the symmetric difference are equinumerous, and
moreover that they are not segregated but intermix in such a manner
that the height function along the border of the modified diamond is
within $o(n)$ of the height function along the border of the true
diamond.  Lastly, let us suppose that the modified diamond has at
least one domino tiling.  Then we can conclude that the average height
function for the modified diamond is within $o(n)$ of the average
height function for the true diamond.

Notice that these results give us no direct information about how
individual placement probabilities change in response to small changes
in the shape of the boundary, though some weak information can be
obtained by way of the height function.  It would be quite interesting
to obtain robustness results for the placement probabilities
themselves.

\subsection{Variance}
\label{ssec-variance}

In \propref{heightformula} of subsection~\ref{ssec-height}~we gave a
formula for the normalized average height function, or rather its
limit as the size $n$ of the Aztec diamond goes to infinity, and in
subsection~\ref{ssec-robustness}~we showed that the same formula
applies to many regions that are roughly similar to the Aztec diamond.
However, we obtained no information about how closely a typical height
function for a region (an Aztec diamond or something else) should
approximate the average height function.  Here we use Azuma's
Inequality \cite[p.~85]{method} to bound the amount of variation that
values of random height functions are likely to exhibit.

Let $H$ denote the (unnormalized) height function corresponding to a
random domino tiling of some simply-connected region in the plane, and
let $v$ denote a vertex in the region, such that there is a path of
$m$ vertices from the boundary of the region to $v$.  We will show in
this subsection that the variance of the random variable $H(v)$ is at
most $64m$.  In fact, we actually get a stronger result:

\begin{thm}
\label{variancethm}
Let $f$ be a partial height function defined on the boundary of a
simply-connected region $R$, and let $v$ be a vertex in the interior
of $R$, such that there is a path of $m$ vertices from the boundary of
$R$ to $v$.  Then, for all $c>0$, the probability that $H(v)$ (the
value of a random height function at $v$ under the uniform
distribution $\mu_f$) differs from its expected value by more than $c
\sqrt{m}$ is less than $2e^{-c^2 / 32}$.
\end{thm}

\begin{proof}
Let $x_0, x_1, \dots, x_{m-1} = v$ be a lattice-path connecting a
point $x_0$ on the boundary of $R$ to the point $v$.  Let $F_k$ be the
partition of the space of possible height functions in which two
height functions are regarded as equivalent if they agree at
$x_0,x_1,\dots,x_{k-1}$.  Let $M_k$ be the conditional expectation
$\Exp[H(v) | F_k]$, the function from the set of height functions to
the reals that assigns to each height function $h$ the average value
of $h'(v)$ as $h'$ ranges over all height functions in the equivalence
class of $h$.

Note that $M_m$ is just the function $H(v)$ itself, while $M_0$ is the
average value of the height at $v$, averaged over all height
functions.  The functions $M_0, M_1, \dots, M_m$ form a martingale;
that is,
\[
\Exp[M_{k+1} | F_k] = M_k.
\]
On each component of $F_k$, $M_{k+1} = \Exp[H(v) | F_{k+1}]$ takes on
at most two distinct values, according to the two different values of
$H(x_k)$ that are consistent with the already-known values of
$H(x_0),H(x_1),\dots,H(x_{k-1})$.  From \corref{concor}, we see that
these two values of $M_{k+1}$ differ by at most 4.  Since $M_k$ is
their weighted average, it follows that $M_{k}$ and $M_{k+1}$ never
differ by more than 4.  Then, applying Azuma's Inequality
(Corollary~2.2 on page~85 of \cite{method}) to the quantities $M_k/4$,
we get
\[
\Prob[|M_m-M_0|/4 > t \sqrt{m} ]  <  2e^{-t^2 / 2}.
\]
Replacing $t$ by $c/4$, we get
\[
\Prob[|H(v)-\Exp[H(v)]| > c \sqrt{m} ]  <  2e^{-c^2 / 32}.
\]
This completes the proof.
\end{proof}

If we are interested in estimating the variance, we can derive a
consequence of the preceding inequality: Assuming for simplicity of
derivation (and without loss of generality) that $\Exp[H(v)]=0$, we
have
\begin{eqnarray*}
\Var[H(v)] & = & \Exp[(H(v))^2] \\
& = & \int_0^\infty \Prob[(H(v))^2 > x] \ dx \\
& = & \int_0^\infty \Prob[|H(v)| > \sqrt{x}] \ dx \\
& < & \int_0^\infty 2e^{-x/(32m)} dx \\
& = & 64m.
\end{eqnarray*}

As a final aside, we mention that our proof of \thmref{variancethm}
also yields the following more general result:
\begin{prop}
\label{large}
Let $R$ be a simply-connected region in the plane and let $v,w$ be
vertices in the interior of $R$, such that there is a path of $m$
vertices from $v$ to $w$, staying entirely within $R$.  Then, for all
$c>0$, the probability that $H(v)-H(w)$ (under the uniform
distribution on domino tilings of $R$) differs from its expected value
by more than $c \sqrt{m}$ is less than $2e^{-c^2 / 32}$.
\end{prop}

\subsection{The arctic circle theorem}
\label{ssec-ACT}

We will now use \thmref{variancethm} to complete the proof of the
arctic circle theorem, which we began in Section~\ref{sec-conclusion}.
We still need to show that the polar regions almost never extend very
far into the interior of the inscribed circle $x^2+y^2=\frac12$
(defined relative to normalized coordinates).  Let $\varepsilon > 0$,
and consider the region $R$ in an Aztec diamond of order $n$ defined
by $x^2+y^2<\frac12-\varepsilon$.

A domino is in the north or south polar region if and only if the
heights on the vertices of the domino are equal to those at the same
locations in the all-horizontal tiling, which is the minimal tiling of
the Aztec diamond (under the partial ordering of tilings induced by
comparison of height functions).  An analogous statement connects the
other two polar regions to the all-vertical tiling, which is the
maximal tiling.  \propref{heightformula} shows that, asymptotically,
the average height function disagrees with the minimal and maximal
height functions within the inscribed circle (although outside of that
circle it agrees with one or the other).  In particular, if a domino
in $R$ is part of the polar regions, then the heights on it differ
{}from the average heights at those locations by an amount at least
proportional to $n$.  (Of course, the constant of proportionality
depends on $\varepsilon$.)  We see that, by taking $c=\sqrt{n}$ in
\thmref{variancethm}, the probability that a domino in $R$ will be
part of the polar regions is exponentially small in $(\sqrt{n})^2=n$.
Since the number of dominos in $R$ is on the order of $n^2$, the
probability that any will be contained in the polar regions is
exponentially small.

We have now proved a slightly stronger version of the arctic circle
theorem than that proved in \cite{circle}.  There, it is shown that
for any $\varepsilon > 0$, for sufficiently large $n$, the boundary of
the temperate zone stays within distance $\varepsilon n$ of the arctic
circle with probability greater than $1-\varepsilon$.  We have shown
that this probability differs from $1$ by an amount exponentially
small in $n$.

\subsection{Heterogeneity}
\label{ssec-heterogeneity}

The arctangent formula gives us an indication of a certain sort of
local homogeneity: places in the tiling that are close together tend
to be governed by the same statistics.  Here we would like to prove a
converse result, and show that within the temperate zone, places in
the tiling that are far apart tend to be governed by different
statistics.  More precisely, we would like to show that within the
temperate zone, the quadruple
$$
(p_n,p_e,p_s,p_w) = ({\mathcal P}(x,y),{\mathcal P}(-y,x),{\mathcal
P}(-x,-y),{\mathcal P}(y,-x))
$$
(whose components are the four placement probabilities near the
location $(x,y)$) uniquely determines $x$ and $y$.  This means, to put
it somewhat fancifully, that if you found yourself stranded somewhere
in the temperate zone of a random domino tiling of a huge checkerboard
colored Aztec diamond, then, provided that you had a compass to tell
you which way was north, you could determine your relative position
within the diamond merely by examining the local statistics of the
tiling.

The heterogeneity claim is not hard to prove, since we know that the
level sets for all four placement probabilities are arcs of ellipses
having a very specific geometry.  In particular, level sets for $p_n$
and $p_s$ are arcs of ellipses that intersect in at most two points,
and these two points have the same $y$-coordinate; similarly, level
sets for $p_e$ and $p_w$ are arcs of ellipses that intersect in at
most two points, and these two points have the same $x$-coordinate.
It follows that if two points have the same probability quadruples,
they must have the same $x$- and $y$-coordinates; that is, the two
points must coincide.

However, we wish to prove more.  Consider that the elements of the
quadruple sum to 1, so that the quadruple has three degrees of
freedom.  However, $x$ and $y$ together embody only two degrees of
freedom, so as $x$ and $y$ sweep through their range of joint allowed
values, the quadruple determined by $x$ and $y$ will not sweep through
the full set of probability vectors of length 4.  On the other hand,
the asymptotic normalized average height function ${\mathcal H}$
introduced in subsection~\ref{ssec-height}~manifests exactly two of
the degrees of freedom of $(p_n,p_e,p_s,p_w)$ in its first-order
derivatives.  What we would like to show is that inside the temperate
zone, the map $(x,y) \mapsto (\frac{\partial {\mathcal H}}{\partial
x}, \frac{\partial {\mathcal H}}{\partial y})$ is one-to-one, and has
as its range the region $\{(s,t): |s|+|t| < 2\}$.  (The possible
significance of this fact will be explained more fully in
Section~\ref{sec-speculations}.)  Putting it differently, we may say
that if one views the graph of the restriction of the function
${\mathcal H}$ to the interior of the temperate zone as a surface,
then the Gauss map from the surface to the sphere is injective.

To prove the claim, we first note that, as discussed in
subsection~\ref{ssec-height}, $\frac{\partial {\mathcal H}}{\partial
x}=2p_w-2p_e$ and $\frac{\partial {\mathcal H}}{\partial
y}=2p_n-2p_s$.  With $y$ fixed and $x$ increasing, $p_e$ increases
while $p_w$ decreases, achieving equality (by symmetry) at $x=0$.
Thus, the sign of $\frac{\partial {\mathcal H}}{\partial x}$ tells us
the sign of $x$, and similarly, the sign of $\frac{\partial {\mathcal
H}}{\partial y}$ tells us the sign of $y$.  Hence, to prove the
injectivity of the map, it suffices to focus on the part of the
temperate zone that lies in the interior of one particular quadrant,
say the second.  Within that quarter-disk, $\frac{\partial {\mathcal
H}}{\partial x}$ and $\frac{\partial {\mathcal H}}{\partial y}$ are
both non-negative functions, taking the values $0$ on the respective
axes $x=0$ and $y=0$ and increasing as one moves away from these axes.
These monotonicity properties do not of themselves rule out the
possibility that $\frac{\partial {\mathcal H}}{\partial x}$ and
$\frac{\partial {\mathcal H}}{\partial y}$ have the same value for two
different points in that quarter-disk, so we must resort to a slightly
more arduous approach.

Using the arctangent formula, one can check that
$$
\frac{\cos\Big(\frac{\pi}{2}\,\frac{\partial{\mathcal H}}{\partial
y}\Big)}
{\cos\Big(\frac{\pi}{2}\,\frac{\partial{\mathcal H}}{\partial x}\Big)}
 = \frac{1-x^2-3y^2}{1-3x^2-y^2}
% The \Big's are used here instead of \left and \right in order to
% ensure that the parentheses in the numerator and denominator have
% the same size.
$$
and
$$
\frac{\sin\Big(\frac{\pi}{2}\,\frac{\partial{\mathcal H}}{\partial
y}\Big)} 
{\sin\Big(\frac{\pi}{2}\,\frac{\partial{\mathcal H}}{\partial x}\Big)}
= -\frac{y}{x}.
% See the note above.
$$
(If $3x^2+y^2=1$, then the first ratio is not defined.  However, since
$x^2+y^2<\frac12$, either the first ratio or its reciprocal is
defined.)  Given the values of these two ratios, there are in general
at most two possibilities for $(x,y)$, only one of which will be in
the desired quadrant.  The only case in which knowledge of the two
ratios does not restrict us to at most two possibilities for $(x,y)$
is when the first ratio is $1$.  This happens iff $\frac{\partial
{\mathcal H}}{\partial x} = \frac{\partial {\mathcal H}}{\partial y}$,
i.e., along the line through the origin that bisects the quadrant.
Since one can check using the explicit formulas for $\frac{\partial
{\mathcal H}}{\partial x}$ and $\frac{\partial {\mathcal H}}{\partial
y}$ that the partial derivatives increase as one moves away from the
axes along that line, they still determine $(x,y)$.  It follows that
the map $(x,y) \mapsto (\frac{\partial {\mathcal H}}{\partial
x},\frac{\partial {\mathcal H}}{\partial y})$ is injective on the
quadrant, as was to be shown.

Now we will see that the map is in fact a surjection to the set
$\{(s,t): |s|+|t|<2\}$.  If one sets $x=(1-t-ct^2)/2$ and
$y=(1+t-ct^2)/2$ with $c > \frac12$ (so that $(x,y)$ lies on a
parabola that is symmetric about the axis $x=y$ and that lies inside
the closed temperate zone in the vicinity of $(\frac12,\frac12)$),
then, sending $t$ to zero from above, we find that the north-going and
east-going probabilities tend towards
\[
\frac{1}{2} + \frac{1}{\pi} \tan^{-1}\frac{1}{\sqrt{2c-1}}
\]
and 
\[
\frac{1}{2} + \frac{1}{\pi} \tan^{-1}\frac{-1}{\sqrt{2c-1}} , 
\]
which sum to 1 for all $c$ between $\frac12$ and infinity and which
vary (as an ordered pair) over the open segment connecting $(1,0)$ to
$(\frac12,\frac12)$, as $c$ goes from $\frac12$ to infinity.  Plugging
the limits $p_n \rightarrow p$, $p_s \rightarrow 0$, $p_e \rightarrow
1-p$, $p_w \rightarrow 0$ into the formulas $\frac{\partial {\mathcal
H}}{\partial x}=2p_w-2p_e$ and $\frac{\partial {\mathcal H}}{\partial
y}=2p_n-2p_s$, we find that the boundary of the open square $\{(s,t):
|s|+|t|<2\}$ consists of limit points of the set of tilts
$(\frac{\partial {\mathcal H}}{\partial x}, \frac{\partial {\mathcal
H}}{\partial y})$ that are achieved by the average height function in
the temperate zone, and hence (by continuity) that that we do indeed
obtain the open square as the set of tilts achieved by the average
height function in the temperate zone.

\subsection{Entropy}
\label{ssec-entropy}

The {\sl entropy} of a random variable that takes on any of $N$ values
with respective probabilities $q_1,\dots,q_N$ is defined as
$\sum_{i=1}^N -q_i \log q_i$ (with $0 \log 0 = 0$ by convention); for
example, the entropy of a uniform random domino tiling of the Aztec
diamond of order $n$ is $\frac{n(n+1)}{2} \log 2$, because there are
exactly $2^{n(n+1)/2}$ tilings (see \cite{alternating} for a proof).
We have seen that for large $n$, nearly all of this entropy is due to
the variety exhibited inside, as opposed to outside, the arctic
circle.  It would be good to have more quantitative information on
this.  Specifically, given a patch of an Aztec diamond, one can define
a random variable whose values are the near-tilings of the patch that
result from restricting a uniform random tiling of the Aztec diamond
to just the patch (such near-tilings are allowed to have untiled
squares along the boundary of the patch), and one can consider the
entropy of this new random variable.  If the patch is very large
(while not long and skinny, for example like a $2 \times n$
rectangle), but the order of the Aztec diamond is much larger still,
then we believe that this entropy, when divided by the area of the
patch, is close to a value which we would call the {\sl local
entropy}, and which would depend only on the normalized location of
the patch.

In this subsection, we make a small start towards calculating local
entropy by showing that it vanishes outside the arctic circle and that
it is positive inside the arctic circle (assuming it is well-defined
there).  Assuming that local entropy is well-defined everywhere, this
gives us another way of interpreting the arctic circle, namely as the
boundary between the zero-entropy region and the positive-entropy
region.

The vanishing (and perforce the well-definedness) of local entropy
outside the temperate zone is a simple consequence of the arctic
circle theorem.  To prove the other half of our claim, consider an $m
\times m$ patch $P$ sitting inside the temperate zone of an extremely
large Aztec diamond, with $m$ even.

If $a$, $b$, $c$, and $d$ are the northwest, northeast, southwest, and
southeast squares in a $2\times 2$ block in a plane region $R$ that
can be tiled by dominos, then the proportion of tilings of $R$ that
have a horizontal domino covering squares $a$ and $b$ and another
horizontal domino covering squares $c$ and $d$ (write this proportion
as $p_{ab,cd}$ for short) is clearly equal to the proportion
$p_{ac,bd}$ of tilings that contain vertical dominos covering squares
$a$ and $c$ and squares $b$ and $d$; moreover, by one of the lemmas
proved in \cite{gip}, both proportions are equal to
$p_{ab}p_{cd}+p_{ac}p_{bd}$, where $p_{ab}$ denotes the proportion of
tilings that have a domino covering $a$ and $b$, etc.\ (that is,
$p_{ab}$, $p_{cd}$, $p_{ac}$, and $p_{bd}$ are just placement
probabilities under uniform random tiling).  In our particular
situation, if one looks inside the patch $P$ taken from the temperate
zone of a large Aztec diamond, all four placement probabilities are
bounded away from zero, say by $\varepsilon>0$, so the probability
that a random tiling contains a $2 \times 2$ block centered at any
particular vertex in $P$ is at least $4\varepsilon^2$.  In particular,
we can look at the $(m/2)^2$ vertices that are at the centers of the
$(m/2)^2$ non-overlapping $2 \times 2$ blocks into which $P$ can be
naturally decomposed.  Using linearity of expectation, we can see that
the expected number of such $2\times 2$ blocks in a random tiling of
the Aztec diamond is at least $m^2 \varepsilon^2$.  However, this
allows us to set a lower bound on the entropy, as measured by the
variety of configurations one sees locally.  For, by freely rotating
these blocks (i.e., changing horizontal blocks to vertical blocks or
vice versa), we can create $2^{m^2 \varepsilon^2}$ other local
patterns, all equally likely.  Standard techniques in information
theory permit one to conclude that the entropy of the near-tiling of
$P$ is at least $\varepsilon^2 \log 2$ times the area of $P$.

\section{Further Results}
\label{sec-further}

Although we have phrased our results in terms of domino tilings, there
is an easy equivalence between domino tilings of finite regions and
dimer configurations on certain finite graphs.  Specifically, if we
replace each square cell by a vertex, and draw an edge connecting any
two vertices whose associated cells are adjacent, then a domino tiling
of a region corresponds to a dimer-cover of the derived graph, that
is, to a set of edges of the derived graph with the property that
every vertex belongs to exactly one of the chosen edges.  In this way,
the study of domino tilings is seen to be equivalent to the study of
dimer-covers, which is one of the better-understood statistical
mechanics models in two dimensions.  Studying domino tilings of
special regions, such as Aztec diamonds, is tantamount to studying the
dimer model in the presence of special boundary conditions.  The
uniformity of the distribution corresponds to a degenerate situation
in which all dimer configurations have the same energy.

There has been surprisingly little work on the behavior of the dimer
model in the presence of general boundary conditions; researchers in
statistical mechanics have tended to study either toroidal (i.e.,
periodic) boundary conditions or boundary conditions that correspond
to domino tilings of a rectangle.  Our work can in a sense be regarded
as a somewhat strange chapter in the study of the dimer model, in
which highly unphysical boundary conditions are imposed.  (Precursors
of this research include \cite{elser}, \cite{gg}, and \cite{remark}.)

In his original article on the dimer model \cite{dimers}, Kasteleyn
considered imposing an energy function that favors one orientation of
dimer over another (horizontal versus vertical).  The authors of
\cite{gip} followed this lead, and showed how their methods also led
to exact results for random domino tilings of the Aztec diamond when
the distribution was skewed towards dominos of a particular
orientation.  Here, we will state the results that follow from
applying the methods of this paper to the case of biased tilings.

Let $p$ be strictly between $0$ and $1$.  For each $n$, there is a
unique probability distribution on the tilings of the Aztec diamond of
order $n$ (in fact, on any simply-connected region) such that given
any tiling of all of the diamond except for a $2 \times 2$ block, the
conditional probability that the $2 \times 2$ block will contain two
horizontal dominos is $p$.  For more details on this distribution, see
\cite{circle} or \cite{gip}.  We call this the Gibbs distribution with
bias $p$.  (For more information on Gibbs distributions in general,
see \cite{georgii}.)

The main difference between the biased distribution and the uniform
distribution is the shape of the temperate zone.  We will see shortly
that, in the biased case, its boundary is given by the ``arctic
ellipse'' $\frac{x^2}{p}+\frac{y^2}{1-p}=1$ (in normalized
coordinates).  It was conjectured in \cite{circle} that the analogue
of the arctic circle theorem holds in the biased case.  Our methods
prove that conjecture, as well as an arctangent formula that describes
the behavior within the temperate zone.

We begin by defining the biased placement probabilities
$\Pl_p(\ell,m;n)$ the same way we defined the ordinary placement
probabilities (except, of course, that we use the biased
distribution).  The biased creation rates are also defined analogously
to the ordinary creation rates, by
$$
\Cr_p(\ell,m;n) = \frac{1}{p}(\Pl_p(\ell,m;n) - \Pl_p(\ell,m-1;n-1)).
$$

The proofs depend on a biased version of \propref{gessel}, which is
proved in \cite{gip}.  To state it, we will need a more general form
of Krawtchouk polynomial.  Define $c_p(a,b;n)$ to be the coefficient
of $z^a$ in $(1+(1-p)z/p)^{n-b}(1-z)^b.$ (See \cite[p.~151]{codes}.)
We need the following result from \cite{gip}:

\begin{prop}
\label{biasgessel}
Let $0 < p < 1$, and set $a = (\ell+m+n)/2$ and $b=(\ell-m+n)/2$.  If
$a,b \in \Z$, then
$$
{\Cr_p}(\ell,m;n+1) =
c_p(a,b;n)c_p(b,a;n)p^n.
$$
Otherwise, ${\Cr_p}(\ell,m;n+1) = 0.$
\end{prop}

Using this proposition, straightforward modifications to the proof of
\propref{crest} prove the following generalization:

\begin{prop}
\label{biascrest}
Fix $\varepsilon > 0$.  If $\frac{\ell^2}{p} + \frac{m^2}{1-p} \le
(1-\varepsilon)n^2$ and $\ell+m \equiv n \pmod{2}$, then
$$
\Cr_p(\ell,m;n+1) = \frac{2\cos^2
\Phi_p(\ell,m;n)}{\pi\sqrt{(p-p^2)n^2-(1-p)\ell^2-pm^2}} +
O_\varepsilon(n^{-2})
$$
for some function $\Phi_p(\ell,m;n)$, which can be determined
explicitly.
\end{prop}

Every result needed for the proof of \thmref{main} (such as the
creation rate estimates outside the arctic ellipse) has a
straightforward generalization to the biased case; in the interest of
saving space, we will omit their statements.  The proofs are
completely analogous to the proofs for the uniform distribution.  One
arrives at the following biased counterpart to \thmref{main}:

\begin{thm}
\label{biasmain}
Let $0 < p < 1$, and let $U$ be an open set containing the points
$(\pm p, 1-p)$.  If $(x,y)$ is the normalized location of a
north-going domino space in the Aztec diamond of order $n$, and $(x,y)
\not\in U$, then, as $n \rightarrow \infty$, the placement probability
at $(x,y)$ for the Gibbs distribution with bias $p$ is within $o(1)$
of ${\mathcal P}_p(x,y)$, where
$$
{\mathcal P}_p(x,y) = 
\begin{cases}
0&\!\!\!\!\hbox{if $\frac{x^2}{p}+\frac{y^2}{1-p}\ge1$ and $y<1-p$,}\\
1&\!\!\!\!\hbox{if $\frac{x^2}{p}+\frac{y^2}{1-p}\ge1$ and $y>1-p$,
and}\\ 
\frac{1}{2}+\frac{1}{\pi}\tan^{-1}\!\!\left(\!\frac{y-(1-p)}
{\sqrt{p-p^2-(1-p)x^2-py^2}}\!\right)\!
&\!\!\!\!\hbox{if $\frac{x^2}{p}+\frac{y^2}{1-p}<1$}.\\
\end{cases}
$$
The $o(1)$ error bound is uniform in $(x,y)$ (for $(x,y) \not\in U$).
\end{thm}

Similarly, the south-going, east-going, and west-going placement
probabilities near $(x,y)$ are approximated by ${\mathcal
P}_p(-x,-y)$, ${\mathcal P}_{1-p}(-y,x)$, and ${\mathcal
P}_{1-p}(y,-x)$, respectively. This follows from \thmref{biasmain} by
rotational symmetry.

One can also prove biased versions of the robustness and variance
results from
subsections~\ref{ssec-robustness}~and~\ref{ssec-variance}.  (In fact,
the proofs are practically identical to the proofs given in those
subsections.)  Using them in combination with the same methods used in
subsection~\ref{ssec-ACT}, we can prove a slightly strengthened
version of the ``arctic ellipse conjecture'' from \cite{circle}:

\begin{thm}
\label{aec}
Let $0 < p < 1$, and $\varepsilon > 0$.  The probability that, in a
random domino tiling with bias $p$ of an Aztec diamond of order $n$,
the boundary of the polar regions is more than a distance
$\varepsilon$ in normalized coordinates from the ellipse
$\frac{x^2}{p}+\frac{y^2}{1-p} = 1$ is exponentially small in $n$.
\end{thm}

\section{Speculations}
\label{sec-speculations}

In this article we have focused primarily on one particular family of
finite regions, namely, Aztec diamonds.  Here we will indicate what it
might mean to have a theory that would apply to {\sl all}
simply-connected finite regions, and how Aztec diamonds might play a
role in the project of classifying the different possible local
behaviors that random tilings of such regions can exhibit away from
their boundaries.

The results of
subsections~\ref{ssec-robustness}~and~\ref{ssec-variance}~tell us that
for any large simply-connected region $R$ that can be tiled by
dominos, height functions associated with random tilings of $R$ will
cluster around their average.  We furthermore know that this average
height function depends in a monotone way on the values of the height
function on the boundary of $R$, and is stable under certain kinds of
slight perturbations of the boundary of $R$.  However, what these
theorems do not tell us is whether this dependence is robust under
scaling as well.  \propref{heightformula} tells us that such
robustness does in fact hold for Aztec diamonds: that is, when one
normalizes two large Aztec diamonds, one finds that the normalized
average height functions are very close to one another.  That this is
true along the boundary is a triviality; that it is true in the
interior is a much subtler property, known to us only as a consequence
of \thmref{main}.

We conjecture that scaling-robustness of height functions is true in
general.  That is, suppose $R_1,R_2,\dots$ are finite,
simply-connected, domino-tileable regions that grow without bound,
such that suitably rescaled copies of the $R_n$'s converge to some
compact subset $R^*$ of the plane.  Moreover, suppose that the height
functions associated with the boundaries of the $R_n$'s, when rescaled
by the same respective amounts, converge to some function on the
boundary of $R^*$.  Then we believe that the average height functions
associated with the $R_n$'s, when rescaled, converge on the interior
of $R^*$ as well as on the boundary to some function ${\mathcal H}$.
If the boundary values behave reasonably (perhaps piecewise smoothness
suffices), then ${\mathcal H}$ should be piecewise smooth (with
reasonably shaped pieces).

Under this picture, we view ${\mathcal H}$ as the solution to a
somewhat strange sort of Dirichlet problem.  We will have more to say
about this analogy shortly, but first we must leave the issues of
large-scale structure (embodied in the average height function) and
discuss the small-scale structure of random tilings.

Consider simply-connected regions $R_1,R_2,\dots$ as above.  In each
region $R_n$, choose a north-going domino space $\sigma_n$ with
normalized location $(x_n,y_n)$ in $R^*$, so that $(x_n,y_n)
\rightarrow (x^*,y^*)$ as $n \rightarrow \infty$, and suppose that the
asymptotic renormalized height function ${\mathcal H}$ is
differentiable at $(x^*,y^*)$.  Assume that $(x^*,y^*)$ is in the
interior of $R^*$ and that ${\mathcal H}$ is ``non-extremal'' at
$(x^*,y^*)$, in the sense that its tilt $(s,t) = (\frac{\partial
{\mathcal H}}{\partial x},\frac{\partial {\mathcal H}}{\partial y})$
satisfies $|s|+|t|<2$.  Then we conjecture that the placement
probabilities at the chosen north-going domino spaces $\sigma_n$
converge.  The arctangent formula tells us that the conjecture is in
fact true for Aztec diamonds.

Note that if we were to replace each $\sigma_n$ by another north-going
domino space $\sigma'_n$ obtained by shifting it by some fixed vector
$(i,j)$ with $i+j$ even, we would get the same point $(x^*,y^*)$ in
the normalized limit.  Hence, the preceding conjecture implies
approximate local translation-invariance for the first-order
statistics governing random tilings of large regions, provided one
stays away from the boundary (and the tilt is non-extremal).

\begin{figure}
\begin{center}
\PSbox{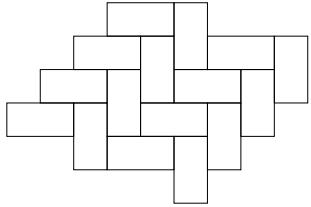 hoffset=-10}{1.5in}{1.5in}
\end{center}
\caption{Herringbone pattern.}
\label{fig-herringbone}
\end{figure}

This corollary gives us a way to understand the importance of our
hypothesis of non-extremality.  For instance, consider the region
shown in Figure~\ref{fig-herringbone}; it has only one tiling, whose
local statistics are in no sense governed by any of the statistics
seen in Aztec diamonds.  Taking a suitable limit of such regions one
gets a height function whose tilt $(s,t)$ satisfies $|s|+|t|=2$ and
hence violates non-extremality.  Indeed, the statistics do not even
exhibit local translation-invariance.  (Note also that for the Aztec
diamond itself, ${\mathcal H}(\cdot,\cdot)$ is extremal at $(x^*,y^*)$
if and only if the asymptotic entropy at normalized location
$(x^*,y^*)$ is zero, which is the case if and only if the asymptotic
density of $2\times 2$ blocks at normalized location $(x^*,y^*)$ is
zero.)

Having made a conjecture about convergence of first-order statistics,
we naturally wonder about higher-order statistics as well.  We
conjecture that in fact all finite-order statistics in the vicinity of
the points $(x_n,y_n)$ stabilize as $n \rightarrow \infty$, yielding
statistics that in some sense ``belong'' to the limit point
$(x^*,y^*)$ (as long as the tilt at $(x^*,y^*)$ is non-extremal).
Then, applying the translation-invariance remark made in the preceding
paragraph, it follows that each $(x^*,y^*)$ determines a process whose
values are domino tilings of the entire plane.  For instance, taking
the $R_n$'s to be Aztec diamonds and the point $(x^*,y^*)$ to be the
center of the normalized diamond, it is natural to conjecture that at
the center of the Aztec diamond of order $n$, the local finite-order
statistics converge to those of the maximal entropy process mentioned
in the Introduction.  (This special case was conjectured in
\cite{circle}.)  Letting $(x^*,y^*)$ vary inside the rescaled
temperate zone, we would get a two-parameter family of tiling-valued
processes; they would all be distinct from one another because they
would have distinct first-order statistics.  The maximal entropy
process would be unique among these processes not only in having the
highest entropy but also in being invariant under the full group of
lattice translations, rather than merely the color-preserving subgroup
of index 2.

It can be shown rigorously that such processes, if they exist, have a
combinatorial analogue of the ``Gibbs property'' studied in
equilibrium statistical mechanics; that is, given a tiling of a
cofinite subset of the plane whose finite complement is tileable, if
one conditions the random process on that particular tiling, then the
conditional distribution on tilings of the entire plane is uniform.

Here we leave aside caution and put forward some conjectures about
what sort of shape the ultimate theory we are striving towards will
take.  These surmises might be false, but we believe they are the
natural avenues to pursue in further investigations of the theory.

In the first place, we conjecture that the tiling-valued processes
associated with the points $(x^*,y^*)$ will turn out to be ergodic, or
indecomposable, in the usual sense of the theory of dynamical systems.
It is not hard to use the ergodic theorem for $\Z^2$-actions (see
\cite{ergodic}) to show that every ergodic, translation-invariant
(under color-preserving translations), tiling-valued random process
determines placement probabilities $p_n$, $p_s$, $p_w$, and $p_e$ and
thence determines a tilt $(s,t)=(2(p_w-p_e),2(p_n-p_s))$.  We predict
that in those cases where the tilt is non-extremal (i.e., $|s|+|t|$ is
strictly less than 2), there is in fact a unique ergodic Gibbs measure
with tilt $(s,t)$.  If this were true, it would have many nice
consequences; for instance, the four numbers $p_n,p_s,p_e,p_w$ would
all be determined by the pair $(s,t)$, and thus would exhibit only two
degrees of freedom, despite the fact that the only obvious constraint
governing them is $p_n+p_s+p_e+p_w=1$.  A further nice property is
that the temperate zones of Aztec diamonds would be universal in the
sense that they would manifest, in the limit, all possible forms of
non-extremal local behavior that random tilings of large
simply-connected regions can manifest away from boundaries.  This
universality is not peculiar to Aztec diamonds, but instead arises
{}from the fact, proved in subsection~\ref{ssec-heterogeneity}, that
Aztec diamonds exhibit all possible non-extremal tilts.

An especially nice benefit of the preceding conjecture is that it
would open the door to a variational approach to the problem of
finding the average height function $\mathcal H$ on $R^*$ given only
its values on the boundary of $R^*$.  Given any candidate for
$\mathcal H$, define $N_n$ as the number of domino tilings of $R_n$
whose normalized height functions stay close to $\mathcal H$.  It does
not seem too far-fetched to hope that the logarithm of $N_n$, when
divided by the area of $R_n$, converges to an integral over $R^*$, in
which the integrand is the entropy associated with the unique ergodic
Gibbs process with tilt $(\frac{\partial {\mathcal H}}{\partial x},
\frac{\partial {\mathcal H}}{\partial y})$.  Since finding the average
height function on $R_n$ corresponds in some sense to maximizing
$N_n$, we would hope that finding the asymptotic normalized height
function on $R^*$ corresponds to maximizing this integral.  It might
not always be possible to solve the associated calculus of variations
problem explicitly, but such a theorem would be a major advance
towards a complete understanding of how the presence of boundary
conditions can affect the behavior of a domino tiling in the interior
of a region.

The preceding idea has in fact been used by physicists, in the context
of crystals; see for example \cite[pp.~3562--3563]{nhb}.  There, it is
claimed that the shape of a crystal surface is determined by
minimizing the total surface free energy, which is obtained by
integrating a local contribution (the surface free energy density)
depending only on the gradient of the surface.  This is believable
physically, but in any particular lattice model it seems difficult to
establish rigorously; it is not even clear on purely mathematical
grounds why there should exist a surface free energy density depending
only on the gradient.  The analogous statement in random tiling theory
is the existence of a local entropy depending only on the tilt of the
height function, but it is conceivable (although we consider it
unlikely) that the local entropy might not be determined by the local
asymptotic behavior of the normalized height function.  The only
approach that we know of that might lead to a rigorous proof (or even
a heuristic argument) is to prove the conjectures above about local
statistics and Gibbs measures.

\section*{Acknowledgements}

We thank Robin Pemantle for suggesting the idea behind
\lemref{monotone}.  Thanks also to David Feldman for providing helpful
comments on the manuscript, to Sameera Iyengar for writing the program
that produced the random tiling shown in Figure~\ref{fig-aztec}, to
M.~Josie Ammer and Dan Ionescu for writing the first programs to
compute placement probabilities for Aztec diamonds, to Pramod Achar,
Federico Ardila, Dan Ionescu, and Ben Raphael for helping to write the
programs that produced the random tilings shown in
Figures~\ref{fig-fake1}, \ref{fig-fake2}, and \ref{fig-real1}, and to
David Wilson for helping to create the figures.

\end{document}